\documentclass[a4paper,11pt]{article}

\usepackage{amsfonts}
\usepackage{epsf,latexsym,amsfonts,amsbsy,mathrsfs}
\usepackage{cite}
\usepackage{amsmath, amssymb, amsthm,amscd,amsxtra}
\usepackage[]{graphicx,subfigure}
\usepackage{float}
\usepackage{multirow}
\usepackage{bbm}
\usepackage{stmaryrd}
\usepackage{epsfig}

\newtheorem{theorem}{Theorem}[section]
\newtheorem{lemma}[theorem]{Lemma}

\theoremstyle{definition}
\newtheorem{definition}[theorem]{Definition}
\newtheorem{example}[theorem]{Example}

\theoremstyle{remark}
\newtheorem{remark}[theorem]{Remark}

\numberwithin{equation}{section}

\newcommand{\eproof}{$\quad \Box$}

\title{Multigrid Methods for Space Fractional Partial Differential Equations\thanks{The work of the first author is supported
by National Natural Science Foundation of China (No.
10901027). The work of second author is supported by the National Basic Research Program under the Grant
2011CB30971 and National Natural Science Foundation of China (No. 11171335, 11225107).}}

\author{Yingjun Jiang\thanks{Department of Mathematics and Scientific
Computing, Changsha University of Science and Technology, Changsha,
410076, China (jiangyingjun@csust.edu.cn).} and Xuejun
Xu\thanks{LSEC, Institute of Computational Mathematics and
Scientific/Engineering Computing, Academy of Mathematics and Systems
Science, Chinese Academy of Sciences, P.O. Box 2719, Beijing,
100190, P.R. China ({\tt xxj@lsec.cc.ac.cn}).}}

\date{}

\begin{document}
\maketitle

\begin{abstract}
We propose some multigrid methods for solving the
algebraic systems resulting from finite element approximations of space
fractional partial differential equations (SFPDEs). It is shown that
our multigrid methods are optimal, which means the convergence
rates of the methods are independent of the mesh size and mesh level. Moreover, our theoretical
analysis and convergence results do not require regularity assumptions of the
model problems.
Numerical results are given to support our theoretical findings.
\end{abstract}

\noindent \hspace{1cm}{\bf Keywords.} fractional differential equations,
multigrid methods, optimal convergence

\section{Introduction}\label{sec1}
 Fractional partial differential
equations (FPDEs) have found many impressive applications  in
lots of fields, such as finance, phase transitions, stratified
materials, anomalous diffusions (see \cite{Nezza1} and references
therein).  To solve them, both analytical and numerical methods
are used in the literature.  The analytical methods like the
Fourier transform method, the Laplace transform method and the
Mellin transform method have been developed to seek closed-form
analytical solutions \cite{Podlubny1}. Since such closed-form
analytical solutions are unavailable in most cases, extensive
researches have already been carried out on the development of
numerical methods for fractional partial differential equations
like finite difference methods (see e.g.,
\cite{Beumer1,Cui1,Langlands1,Meerschaert1,Meerschaert2,Sousa1,Tadjeran1}),
finite element methods (see e.g., \cite{Deng1,Ervin0,Liu1}), and
spectral methods \cite{Li1,Lin2}.

Let $\Omega$ be a polyhedral
domain in $\mathbb{R}^d$, we consider the space fractional
partial differential equations (SFPDEs):
find $u(x)$ such
that (see \cite{Ervin2})
\begin{eqnarray}\label{eqn source}
-\int_{S^{d-1}}D^{2\alpha}_{z}u(x)\, \tilde{M}(z)dz+cu(x)=f(x), \quad x\in
\Omega,
\end{eqnarray}
\begin{equation}\label{eqn boundary}
u|_{R^d\setminus \Omega}=0,
\end{equation}
where $1/2< \alpha \leq 1$, $c\geq 0$, $f$  is a source term,
$||\cdot||_2$ denotes the standard  Euclidean norm, $S^{d-1}=\{z\in
R^d; ||z||_2=1\}$, $\tilde{M}(z)$ is a probability density function on
$S^{d-1}$, and $D^{2\alpha}_{z}$, which will be given later, denotes
the directional derivative of order $2\alpha$ in the direction of
the unit vector $z$. Here we assume $\tilde{M}$ is symmetric about origin,
i.e., $\tilde{M}(z)=\tilde{M}(z')$ if $z,z'\in S^{d-1}$ satisfy $z+z'=0$, which
means that
 the considered problem is a symmetric one.

 One special case of (\ref{eqn source}) is
\begin{equation}\label{eqn dis}
  -\sum\limits_{i=1}^d(p_i\ _{-\infty}D_{x_i}^{2\alpha}+q_i\ _{x_i}D_\infty^{2\alpha})u+cu=f
\end{equation}
and $p_i,q_i\geq0$ satisfying $p_i=q_i$ and $\sum_{i=1}^d(p_i+q_i)=1$, where $_{-\infty}D_{x_i}^{2\alpha},
{_{x_i}D_\infty^{2\alpha}}$ denote Riemann-Liouville fractional
derivatives. Actually, (\ref{eqn dis}) can be obtained from (\ref{eqn
source}) by taking  $\tilde{M}=\sum_{i=1}^d
p_i\delta(z-e_i)+q_i\delta(z+e_i)$, where $e_i$ is the $i$th
column of identity matrix in $\mathbb{R}^{d\times d}$ and $\delta$
the Dirac function on $S^{d-1}$.
The corresponding time-dependent equation of (\ref{eqn source}) can be used to describe a general super-diffusion process (see
\cite{Meerschaert0}), which is
 an appropriate extension from one dimensional problem
\begin{equation}\label{eqn surce one d}
   \frac{\partial u}{\partial t}-(p\ _{-\infty}D_x^{2\alpha}+q\ _xD_\infty^{2\alpha})u+cu=f.
\end{equation}
 As to the super-diffusion, please refer
to \cite{Metzler1} for details.

One of the greatest challenges for numerically solving SFPDEs is how
to reduce the computation costs.  Due to the nonlocal properties of
fractional differential operators, numerical methods for linear
SFPDEs tend to yield the linear
equations $Ax=b$ with the following characteristics: 1). the
coefficient matrix $A$ is dense  or full; 2). the condition number
of $A$ increases fast, as the mesh becomes fine. Reducing the
computation costs for SFPDEs is harder than doing it for the integer
order PDEs. Some methods have already been designed to overcome this
difficulty, such as alternating-direction implicit methods
(ADI)\cite{Meerschaert2,Wang2,Wang3}, and iterative methods \cite{Lin1,Pang1,Wang3,Wang4,Wang5,Wang6}.

Iterative methods seem to be efficient tools for solving SFPDEs. Actually
two issues in this situation need to be concerned for efficiency:  one is to do the matrix-vector
multiplications efficiently, and the other is to find good
preconditioners.  As to the first issue, some literatures are contributed: in \cite{Wang1}, with the notice of
 Toeplitz-like structure of the coefficient matrix, the matrix-vector multiplications are done with $O(N\log
N)$ complexity by using a fast Fourier
transform (FFT) \cite{Chan1,Chan2}.  This technique of "matrix-vector multiplication" has
been widely used to improve the
efficiency of iterative methods for the SFPDEs
\cite{Lin1,Pang1,Wang3,Wang4,Wang5,Wang6}.
As regards the second issue, some literatures should be listed as follows: the first relevant paper may be \cite{Bramble1} in which   a multilevel preconditioner of fractional power was put forward; in \cite{Lin1}, the authors propose preconditioners constructed
by some banded matrices of fixed band width; in \cite{Wang5},
the authors present a preconditioner by some symmetric positive Toeplitz matrixs;  moreover a new preconditioner is designed in \cite{Lei1} through some  circulant matrixs.

It is known that multigrid methods are optimal iterative procedures, which have been widely used for integer order
PDEs (see e.g., \cite{Brammble,Smith1}).  In recent years, some researchers begin to investigate
multigrid methods for solving SFPDEs. For instance, in  \cite{Zhou1}, Zhou and Wu  apply the
multigrid method  to solve one dimensional steady SFPDEs, and
in \cite{Pang1}, the authors consider the V-cycle multigrid method for solving
corresponding time-dependent problems. But till now, no satisfactory
convergence results have been obtained for the multigrid methods for
solving SFPDEs. Actually, in \cite{Pang1}, the authors only conduct the theoretical analysis
for the two-level multi-grid method,  and Zhou and Wu in \cite{Zhou1} get the
convergence results only under the assumption that the adjoint problem hold
sufficiently smooth solution.

In this paper, we introduce a V-cycle multigrid method with one smoothing step on each level  to solve
linear algebraic systems resulting from the finite element approximations of the SFPDEs
(\ref{eqn source}). It is shown that
our V-cycle multigrid methods are optimal, which means the convergence
rates are independent of the mesh size and mesh level. Moreover, our theoretical
analysis and the convergence results in this paper do not require any regularity assumptions of the
model problems. To the best of our knowledge, this paper is a first attempt to give a rigorous theoretical
analysis for the V-cycle multigrid methods for the finite element approximations of SFPDEs in any dimensions.

This paper is also the first work to design the fast solver for the SFPDE (\ref{eqn source}) with $M$ being
 a continuous function. Among the current numerical methods for SFPDEs, most of them are for one
dimensional problems and for some special high dimensional problems like
(\ref{eqn dis}), and only a few are for more general problems like (\ref{eqn
source}).
Actually, only \cite{Ervin2,Roop1}
study the numerical methods for (\ref{eqn source}):   in
\cite{Ervin2}, the authors consider the finite element approximation
for (\ref{eqn source}) and in \cite{Roop1}, the author studies the corresponding time-dependent case.

 In the rest of the paper, no loss of generality, we restrict ourselves to the case $d=2$, namely,
 we consider the problem (\ref{eqn source}) in $\mathbb{R}^2$. For $\Lambda\subset \mathbb{R}^2$,
 denote $L^2(\Lambda)$ the space of all measurable function $v$ on $\Lambda$ satisfying
$\int_\Lambda (v(x))^2dx<\infty,$ and $C^\infty_0(\Lambda)$  the
space of infinitely differentiable functions with compact support in
$\Lambda$. Set
$$(v,w)_\Lambda=\int_\Lambda vw dxdy,\quad ||v||_\Lambda=(v,v)^{1/2}_\Lambda,$$
and they are abbreviated as $(v,w)$ and $||v||$ respectively if
$\Lambda=\mathbb{R}^2$.

To simplify our statement, we make a convention here: function $v$
defined on a domain $\Lambda\subset\mathbb{R}^2$ also denotes its
extension on $\mathbb{R}^2$ which extends $v$ by zero outside
$\Lambda$. The constant $C$  with or without subscript will denote
a generic positive constant which may take on different values in
different places. These constants will always be independent of  the mesh
sizes and  levels in the multigrid methods.
Following \cite{Xu1},  we also use symbols $\lesssim,\gtrsim $ and
$\approx $ in this paper. That $a_1\lesssim b_1$, $a_2\gtrsim b_2$
and $a_3\approx b_3$ mean that $a_1\leq C_1 b_1$, $a_2\geq C_2 b_2$
and $C_3b_3\leq a_3\leq C'_3 b_3$ for some positives $C_1,C_2,C_3$
and $C'_3$.

The rest of the paper is organized as follows: for the sake of
completeness, in section \ref{sec}, we give our model problem
 and the corresponding finite element
discretization. In section \ref{sec3}, we present our V-cycle multigrid methods and introduce some basic
theoretical results.  In section \ref{sec4}, we shall prove
the convergence of the multigrid methods. Finally in section
5, the numerical results are given to verify our theoretical
findings.

\section{The model problem and its discretization}
In this section, we shall present the SFPDE  in $\mathbb{R}^2$, and
then introduce its variational formulation and corresponding finite element
discretization.

\subsection{The model problem} \label{sec}
 We first introduce the concepts of directional integrals and
derivatives \cite{Ervin2}.
\begin{definition}\label{defi 1}\cite{Ervin2}
Let $\mu>0$, $\theta\in \mathbb{R}$. The $\mu$th order fractional
integral in the direction $z=(\cos\theta,\sin\theta)$ is defined by
$$D^{-\mu}_z v(x,y):=D^{-\mu}_\theta v(x,y)=\int_0^\infty \frac{\tau^{\mu-1}}{\Gamma(\mu)}v(x-\tau \cos\theta,y-\tau\sin\theta)d\tau,$$
where $\Gamma$ is the Gamma function.
\end{definition}

\begin{definition}\label{defi 2}\cite{Ervin2}
Let $n$ be a positive integer, and $\theta\in \mathbb{R}$. The $n$th
order derivative in the direction of $z=(\cos\theta,\sin\theta)$ is
given by
$$D^n_\theta v(x,y):=\left(\cos\theta\frac{\partial}{\partial x}+\sin\theta\frac{\partial}{\partial y}\right)^nv(x,y).$$
\end{definition}
\begin{definition}\label{defi 3}\cite{Ervin2}
Let $\mu>0$, $\theta\in \mathbb{R}$. Let $n$ be the integer such
that $n-1\leq\mu<n$, and define $\sigma=n-\mu$. Then the $\mu$th
order directional derivative in the direction of
$z=(\cos\theta,\sin\theta)$ is defined by
$$D^\mu_z v(x,y):=D^\mu_\theta v(x,y)=D^n_\theta D_\theta^{-\sigma}v(x,y).$$
\end{definition}
If  $v$ is viewed as a function in $x$,
$D^{\mu}_0$, $D^{\mu}_\pi$ are
just the left and the right Riemamm-Liouville derivatives (see e.g., \cite{Podlubny1,Samko1}). The
fractional derivative operators in problem (\ref{eqn source}) are
related to the following  fractional
derivative:
\begin{definition}\label{defi 4}\cite{Ervin2}
Assume that $v:\mathbb{R}^2\rightarrow \mathbb{R}$, $\mu>0$. The $\mu$th
order fractional derivative with respect to the measure $\tilde{M}$ is
defined as
$$D^\mu_{\tilde{M}}v(x,y):=\int_{S^1}D^\mu_\theta v(x,y)\tilde{M}(\theta)d\theta,$$
where $S^1=[0+\nu, 2\pi+\nu)$ with a suitable scalar $\nu$, and
$\tilde{M}(\theta)$, which satisfies
$\int_{\nu}^{2\pi+\nu}\tilde{M}(\theta)d\theta=1$, is a periodic function
with period $2\pi$. Usually we take $\nu=0$, if it causes no
unreasonable expression (see (\ref{eqn a0})).
\end{definition}
\begin{remark}\label{remark 1}
It is easy to check that
$$D^2_{\tilde{M}} v(x,y)=a_{11}\frac{\partial^2 v}{\partial x^2}+a_{22}\frac{\partial^2 v}{\partial y^2}+2a_{12}\frac{\partial^2 v}{\partial x\partial y},$$
where $a_{11}=\int_{0}^{2\pi}\cos^2\theta \tilde{M}(\theta)d\theta$,
$a_{22}=\int_{0}^{2\pi}\sin^2\theta  \tilde{M}(\theta)d\theta$ and
$a_{12}=2\int_{0}^{2\pi}\cos\theta\sin\theta  \tilde{M}(\theta)d\theta$ (see
also \cite{Meerschaert00}). Denote $L$ a positive integer, let
$\theta_k\in [0,2\pi)$ and $p_k\geq 0$, $k=1,2,\ldots,L$, satisfy
$\sum_{k=1}^Lp_k=1$. Assume that $D^{\mu}_\theta v$ is continuous in
$\theta$, and then
\begin{equation}\label{eqn D_M 1}
D^\mu_{\tilde{M}} v=\sum_{k=1}^Lp_kD^\mu_{\theta_k} v(x,y),
\end{equation}
if \begin{equation}\label{eqn a0}
\tilde{M}=\sum_{k=1}^Lp_k\delta(\theta-\theta_k),
\end{equation} where $\delta$ denotes Dirac delta function.
\end{remark}

For $u:\mathbb{R}^2\rightarrow \mathbb{R}$, define differential
operator $L_\alpha$ in $\mathbb{R}^2$ as
$$L_\alpha u=-D^{2\alpha}_{\tilde{M}} u+cu.$$
Denote $\Omega$ a polygonal domain in $\mathbb{R}^2$, set
$1/2<\alpha\leq 1$, and then the model problem of this paper
is to find $u:\bar\Omega\rightarrow \mathbb{R}$ such that
\begin{equation}\label{eqn source new}
\left\{\begin{array}{cc}L_\alpha u=f, & in\ \Omega,\\
 u=0, &  \hbox{on}\
\partial\Omega,\end{array}\right.
\end{equation}
where $f$ is a source term and we assume that $\tilde{M}(\theta)$ satisfies
$\tilde{M}(\theta)=\tilde{M}(\theta+\pi)$ for $\theta\in \mathbb{R}$, i.e.,
(\ref{eqn source new}) is a symmetric problem. Here, we recall the
convection made in Section \ref{sec1}, i.e.,  $u$ also denotes its
extension by zero outside $\Omega$.

\subsection{The variational formulation}
\begin{definition}\label{defi 8}\cite{Tartar1} Let $\mu\geq 0$, $\mathcal{F}{v}(\xi_1,\xi_2)$ be the Fourier transform of $v(x,y)$, $|\xi|=\sqrt{\xi_1^2+\xi_2^2}$.
Define norm
$$||v||_{H^\mu(\mathbb{R}^2)}:=\left\|(1+|\xi|^2)^{\mu/2}|\mathcal{F}{v}|\right\|.$$
Let $H^\mu(\mathbb{R}^2):=\{v\in L^2(\mathbb{R}^2);
||v||_{H^\mu(\mathbb{R}^2)}< \infty\}$.
\end{definition}
For $v\in H_0^\mu(\Omega)$,  we also denote  $||v||_{H^\mu(\mathbb{R}^2)}$ by $||v||_{H^\mu(\Omega)}$. It is known that $H^\mu(\mathbb{R}^2)$ is a Hilbert space equipped
with the inner product
$(v,w)_{H^\mu(\mathbb{R}^2)}=((1+|\xi|^2)^{\mu}\mathcal{F}{v},\overline{\mathcal{F}{w}})$
and $C_0^\infty(\mathbb{R}^2)$ is dense in $H^\mu(\mathbb{R}^2)$
(see \cite{Tartar1}).
Now, we introduce and prove some useful results for the fractional
directional derivatives of functions in $C_0^\infty(\mathbb{R}^2)$.
\begin{lemma}\label{lemma of fourier transform11}\cite{Ervin2}
For $\mu\in \mathbb{R}$, $v\in C^\infty_0(\mathbb{R}^2)$, the Fourier transform of
$D^{\mu}_\theta v$ is
$$\mathcal{F}(D^{\mu}_\theta v(x,y))=\left(2\pi i(\xi_1\cos\theta+\xi_2\sin\theta)\right)^{\mu}\mathcal{F}{v}(\xi_1,\xi_2).$$
\end{lemma}
\begin{lemma}\label{lemma a5}
For $\mu,s>0$, $v,w\in C_0^\infty(\mathbb{R}^2)$,
$$(D^\mu_\theta v,w)=(D^{\mu-s}_\theta v,D^s_{\theta+\pi} w),$$
where $D^{0}_\theta v=v.$
\end{lemma}
\proof By lemma \ref{lemma of fourier transform11} and (\ref{eqn aaaa6}), we know $\mathcal{F}D^{\mu}_\theta v=\left(2\pi i(\xi_1\cos\theta+\xi_2\sin\theta)\right)^{\mu}\mathcal{F}v$, $\mathcal{F}D^{\mu-s}_\theta v=\left(2\pi i(\xi_1\cos\theta+\xi_2\sin\theta)\right)^{\mu-s}\mathcal{F}v$, $\mathcal{F}D^{s}_{\theta+\pi} w=\overline{\left(2\pi i(\xi_1\cos\theta+\xi_2\sin\theta)\right)^{s}}\mathcal{F}w$. Then the lemma follows by Parseval's formula.
\eproof

We define the weak fractional directional
derivative according to the relation $(D^\mu_\theta
v,w)=(v,D^\mu_{\theta+\pi} w)$ which is a special case of Lemma
\ref{lemma a5} (see also Lemma 5.7 in \cite{Ervin2}). Let
$L_{loc}^1(\mathbb{R}^2)$ denote the set of locally integrable
functions on $\mathbb{R}^2$.
\begin{definition}\label{defi 0} Given $\mu>0$,
$\theta\in \mathbb{R}$, let $v\in L^2(\mathbb{R}^2)$. If there is a
function $v_\mu\in L^1_{loc}(\mathbb{R}^2)$ such that
$$(v,{D^\mu_{\theta+\pi}}w)=({v_\mu},w),\quad \forall w\in C_0^\infty(\mathbb{R}^2),$$
then $v_\mu$ is called the weak $\mu$th order derivative in the
direction of $\theta$ for $v$, denoted by ${D^\mu_{\theta}}v$, i.e.,
$v_\mu={D^\mu_{\theta}}v$.
\end{definition}

It is not hard to see that the weak derivative ${D^\mu_{\theta}}v$
is unique if it exists and that the weak derivative coincides with
the correspondent derivative defined in Definition \ref{defi 3} if $v\in
C^\infty_0(\mathbb{R}^2)$.  In the following, we use
${D^\mu_{\theta}}v$ to denote the weak derivative.

\begin{lemma}\label{lemma a9}
Let $\mu>0$. For any $v\in H^\mu(\mathbb{R}^2)$, $0<s\leq\mu$ and $\theta\in \mathbb{R}$, the weak derivative
${D^s_{\theta}}v$ exists and satisfies
\begin{equation}\label{r1}\mathcal{F}D^s_\theta v(\xi_1,\xi_2)=(2\pi i\xi_1\cos \theta+2\pi i\xi_2\sin\theta)^s \mathcal{F}v(\xi_1,\xi_2),\end{equation}
\begin{equation}\label{r2}||D^s_\theta v||\leq C||v||_{H^\mu(\mathbb{R}^2)}.\end{equation}
\end{lemma}
\proof Since $C_0^\infty(\mathbb{R}^2)$ is dense in $H^\mu(\mathbb{R}^2)$, there is a Cauchy sequence $\{v_n\}\subset C_0^\infty(\mathbb{R}^2)$ such that $||v_n-v||_{H^\mu(\mathbb{R}^2)}\rightarrow 0$ as $n\rightarrow 0$.
By lemma \ref{lemma of fourier transform11}, $\mathcal{F}D^{s}_\theta w=\left(2\pi i(\xi_1\cos\theta+\xi_2\sin\theta)\right)^{s} \mathcal{F}{w}$ for $w\in C_0^\infty(\mathbb{R}^2)$. By Parseval's formula and $0<s\leq\mu$, it is not hard to see that $||D^{s}_\theta w||=||\mathcal{F}D^{s}_\theta w||\leq C||w||_{H^\mu(\mathbb{R}^2)}$. So we have $||D^{s}_\theta v_n-D^{s}_\theta v_m||\leq C||v_n-v_m||_{H^\mu(\mathbb{R}^2)}$ and $\{D^{s}_\theta v_n\}$ is a Cauchy sequence in $L^2(\mathbb{R}^2)$.
Denote $v_s\in L^2(\mathbb{R}^2)$ the function to which $\{D^{s}_\theta v_n\}$ converges to. By Lemma
\ref{lemma a5}, for any $w\in C_0^\infty(\mathbb{R}^2)$,
$$(v_n,D^s_{\theta+\pi} w)=(D^{s}_\theta
v_n, w).$$
Take the limits of both sides of the above equation, we obtain
$(v,D^s_{\theta+\pi} w)=(v_s, w)$ for any $w\in C_0^\infty(\mathbb{R}^2)$. So ${D^s_{\theta}}v$ exists and is equal to $v_s$ by Definition \ref{defi 0}. By the definition of Fourier transform for the function in $L^2(\mathbb{R}^2)$,
\begin{equation}\label{r3}
(\left(2\pi i(\xi_1\cos\theta+\xi_2\sin\theta)\right)^{s} \mathcal{F}v_n,v)=(D^s_\theta v_n,\mathcal{F}v), \quad\forall v\in C^\infty_0(\mathbb{R}^2).
\end{equation}
 Because
 $$||v_n-v||_{H^\mu(\mathbb{R}^2)}=\left\|(1+|\xi|^2)^{\mu/2}|\mathcal{F}{(v_n-v)}|\right\|\rightarrow 0,$$
 it is not hard to see that  $\left(2\pi i(\xi_1\cos\theta+ \xi_2\sin\theta)\right)^{s} \mathcal{F}v_n$ converges to
$\left(2\pi i(\xi_1\cos\theta+\xi_2\sin\theta)\right)^{s} \mathcal{F}v$ in $L^2(\mathbb{R}^2)$.
Take the limits of both sides of (\ref{r3}), we obtain (\ref{r1}) by the definition of Fourier transform. (\ref{r2}) can be proved directly by (\ref{r1}) and Parseval's formula. \eproof

\begin{lemma}\label{lemma a7}
Let $\mu,s>0$ with $\mu-s>0$. For $v,w\in H^{\mu+s}(\mathbb{R}^2)$,
\begin{equation}\label{eqn aa4}
(D^\mu_\theta v,D^\mu_{\theta+\pi} w)=(D^{\mu+s}_\theta
v,D^{\mu-s}_{\theta+\pi} w).
\end{equation}
\end{lemma}
\proof For any $g\in H^{\mu+s}(\mathbb{R}^2)$, $||D^\mu_\theta g||$, $||D^\mu_{\theta+\pi} g||$, $||D^{\mu+s}_\theta g||$ and $||D^{\mu-s}_{\theta+\pi} g||$ are all bounded by $C||g||_{H^\mu(\mathbb{R}^2)}$ by Lemma \ref{lemma a9}. Then the lemma follows from that $C_0^\infty(\mathbb{R}^2)$ is dense in $H^{\mu+s}(\mathbb{R}^2)$ and Lemma \ref{lemma a5}.\eproof

Assume that the solution $u$ of (\ref{eqn source new}) is sufficiently
smooth (indeed, that $u\in C^2(\Omega)$ with $u|_{\partial\Omega}=0$
is sufficient). Multiplying both sides of the first equation in (\ref{eqn source new})
with $v\in C^\infty_0(\Omega)$  and integrating over $\Omega$ give
\begin{equation}\label{eqn aaaaaa_1}
-\int_0^{2\pi}(D^{2\alpha}_\theta u,v)\tilde{M}(\theta)d\theta+c(u,v)=(f,v),
\quad v\in C^\infty_0(\Omega).
\end{equation}
Then employing the relation $(D^1_\theta w,v)=(w,D^1_{\theta+\pi}
v)$ (it can be obtained by integration by parts), we obtain
\begin{equation}\label{eqn aaaaaa0}
-\int_0^{2\pi}(D^{2\alpha-1}_\theta u,D^1_{\theta+\pi}
v)\tilde{M}(\theta)d\theta+c(u,v)=(f,v), \quad v\in C^\infty_0(\Omega).
\end{equation}
Then by Lemma \ref{lemma a7}, (\ref{eqn aaaaaa0}) can be rewritten
as
\begin{equation}\label{eqn aaaaaa1}
-\int_0^{2\pi}(D^{\alpha}_\theta u,D^{\alpha}_{\theta+\pi}
v)\tilde{M}(\theta)d\theta+c(u,v)=(f,v), \quad v\in C^\infty_0(\Omega).
\end{equation}

Define the bilinear form $\tilde{B}: H^\alpha_0(\Omega)\times
H^\alpha_0(\Omega)\rightarrow \mathbb{R}$ as
$$
\tilde{B}(u,v):=-\int_0^{2\pi}(D^\alpha_\theta u,D^\alpha_{\theta+\pi}
v)\tilde{M}(\theta)d\theta+c(u,v).
$$
By $\tilde{M}(\theta)=\tilde{M}(\theta+\pi)$ for $\theta\in \mathbb{R}$, it is easy
to check that $\tilde{B}(v,w)$ is a symmetric bilinear  form, i.e.,
$\tilde{B}(v,w)=\tilde{B}(w,v)$ for $v,w\in H^{\alpha}_0(\Omega)$.
The variational formulation of (\ref{eqn source new}) is (see also
\cite{Ervin2}) to find $u\in H^{\alpha}_0(\Omega)$ such that
\begin{equation}\label{eqn variational form new}
\tilde{B}(u,v)=(f,v),\quad \forall v\in H^{\alpha}_0(\Omega).
\end{equation}
Now we restate some results in \cite{Ervin2} about the solvability of (\ref{eqn variational form new}). To guarantee the existence of the solution of (\ref{eqn variational
form new}), we assume that  $\tilde{M}(\theta)$
satisfies
\begin{equation}\label{eqn a1}
\int_0^{2\pi}|(\xi_1\cos\theta+\xi_2\sin\theta)|^{2\alpha}\tilde{M}(\theta)d\theta\geq C_0|\xi|^{2\alpha}
\end{equation}
 for some positive $C_0$. Denote $\kappa=2\pi (\xi_1\cos\theta+\xi_2\sin\theta)$, $E_1=\{(\xi_1,\xi_2)\in \mathbb{R}^2:\xi_1\cos\theta+\xi_2\sin\theta>0\}$, $E_2=\{(\xi_1,\xi_2)\in \mathbb{R}^2:\xi_1\cos\theta+\xi_2\sin\theta<0\}$,
and then by Parseval's formula and Lemma \ref{lemma a9},
\begin{eqnarray}\label{eqn r1}
(D^\alpha_\theta v,D^\alpha_{\theta+\pi}
v)&=&((i\kappa)^{2\alpha} \mathcal{F}v,\overline{\mathcal{F}v})\nonumber\\
&=&(|\kappa|^{2\alpha} \exp({i\alpha\hbox{sign}(\kappa)\pi}) \mathcal{F}v,\overline{\mathcal{F}v})\nonumber\\
&=&(|\kappa|^{2\alpha} \exp({i\alpha\pi}) \mathcal{F}v,\overline{\mathcal{F}v})_{E_1}+(|\kappa|^{2\alpha} \exp({-i\alpha\pi}) \mathcal{F}v,\overline{\mathcal{F}v})_{E_2}\nonumber\\
&=&\cos(\alpha\pi)(|\kappa|^{2\alpha}\mathcal{F}v,\overline{\mathcal{F}v})\nonumber\\
&&\quad +i\sin(\alpha\pi)\left((|\kappa|^{2\alpha}\mathcal{F}v,\overline{\mathcal{F}v})_{E_1}-(|\kappa|^{2\alpha}\mathcal{F}v,\overline{\mathcal{F}v})_{E_2}
\right)\nonumber\\
&=&\cos(\alpha\pi)(|\kappa|^{2\alpha}\mathcal{F}v,\overline{\mathcal{F}v}),
\end{eqnarray}
where the computation of complex please refer to Appendix, in the fourth equality, the Euler formula $\exp(i\kappa)=\cos(\kappa)+i\sin(\kappa)$ is used, the last equality is because the value of $(D^\alpha_\theta v,D^\alpha_{\theta+\pi}v)$ is real and the imaginary part must be zero (another proof for this equality please refer to \cite{Ervin2}). Furthermore, by (\ref{eqn a1}) and $\cos(\alpha\pi)<0$
\begin{eqnarray}\label{eqn r2}
&&-\int_0^{2\pi}(D^\alpha_\theta v,D^\alpha_{\theta+\pi}
v)\tilde{M}(\theta)d\theta\nonumber\\
&=&-\cos(\alpha\pi)\iint_{\mathbb{R}^2}|\mathcal{F}v|^2\int_0^{2\pi}|2\pi (\xi_1\cos\theta+\xi_2\sin\theta)|^{2\alpha}\tilde{M}(\theta)d\theta\, d\xi_1d\xi_2\nonumber\\
&\gtrsim&\iint_{\mathbb{R}^2} |\xi|^{2\alpha}|\mathcal{F}v|^2 d\xi_1d\xi_2.
\end{eqnarray}
For $v\in H^{\alpha}_0(\Omega)$, we have
\begin{eqnarray}\label{r4}||v||^2\leq C_1||D^\alpha_\theta v||^2&=&C_1\iint_{\mathbb{R}^2}\left|2\pi (\xi_1\cos\theta+\xi_2\sin\theta)\right|^{2\alpha} |\mathcal{F}v|^2 d\xi_1d\xi_2\nonumber\\
&\leq& C_2\iint_{\mathbb{R}^2} |\xi|^{2\alpha}|\mathcal{F}v|^2 d\xi_1d\xi_2,
\end{eqnarray}
where the inequality is by (5.15) in \cite{Ervin2} and the equality is by Parseval's formula.
With the combination of (\ref{eqn r2}) and (\ref{r4}), we conclude under condition (\ref{eqn a1}),
\begin{equation}\label{r5}
\tilde{B}(v,v)\gtrsim ||v||^2_{H^{\alpha}(\Omega)},\quad v\in H^{\alpha}_0(\Omega).
\end{equation}
By Lemma \ref{lemma a9}, it is easy to verify that
\begin{equation}\label{r6}
\tilde{B}(v,w)\lesssim ||v||_{H^{\alpha}(\Omega)}||w||_{H^{\alpha}(\Omega)},\quad v,w\in H^{\alpha}_0(\Omega).
\end{equation}
By (\ref{r5}) and (\ref{r6}), using Lax-Milgram theorem, we know that the variational
formulation (\ref{eqn variational form new}) admits a unique solution in
$H^\alpha_0(\Omega)$.

\begin{remark}\label{remark 3}
Condition (\ref{eqn a1}) is easily satisfied. For example, it holds if
$\tilde{M}(\theta)$ is non-zero over a connected set of positive measure in
$[0,2\pi)$ (see \cite{Ervin2}), and it holds when
$\tilde{M}(\theta)=\sum_{k=1}^4p_k\delta(\theta-k\pi/2)d\theta$, with
$p_k\geq 0$ and $p_1+p_3=1$, $p_2+p_4=1$.
\end{remark}

\subsection{The finite element discretization}
  Let ${\cal T}_h$ be a quasi-uniform triangulation of $\Omega$ such that $\bar\Omega=\cup_{K\in {\cal T}_h} K$,
 $h_K$ be the maximal length of the sides of the triangle $K$ and $h=\max_{K\in {\cal T}_h}h_K$.
Denote $P_l(K)$, $l\geq1$, the space of polynomials of degree less
than or equal to $l$ on $K\in {\cal T}_h$. Define the finite dimensional
subspace $V$ associated with ${\cal T}_h$ as
$$
V:=\{v\in  C^0(\bar\Omega): v|_{\partial\Omega}=0, v|_K\in
P_{l}(K),\forall K\in {\cal T}_h\}.
$$
 It is known that
$V\subset H^1_0(\Omega)\subset
H^\alpha_0(\Omega)$. Thus the finite element approximation for (\ref{eqn
variational form new}) is to find $\tilde{u}_h\in V$ such that
 \begin{equation}\label{eqn finite element new r1}
 \tilde{B}(\tilde{u}_h,v)=(f,v),\quad \forall v\in V.
 \end{equation}
The error estimates for the finite element solution  $\tilde{u}_h$ are given in \cite{Ervin2}.

In practical applications, we use the finite element discretization (\ref{eqn finite element new r1}) only when the probability density function $\tilde{M}$ has the discrete form as that in (\ref{eqn a0}) (when $\tilde{M}(\theta)$ is a
continuous function, the finite element discretization (\ref{eqn finite element new r1}) can hardly be realized). For the case that $\tilde{M}(\theta)$ is the continuous function, we propose an alternative finite element discretization instead of (\ref{eqn finite element new r1}).
Here we focus on the case $\tilde{M}(\theta)\in C^1[0,2\pi]$ is a
periodic function with period $2\pi$ to present our alternative  finite element problem:
find $\bar{u}_h\in V$ such that
 \begin{equation}\label{eqn finite element new r2}
 \bar{B}(\bar{u}_h,v)=(f,v),\quad \forall v\in V,
 \end{equation}
where $\bar{B}(\cdot,\cdot)$ is an approximation of $\tilde{B}(\cdot,\cdot)$.
Exactly in this paper,  set a positive integer $N_\theta$ such that
 $N_\theta$ is  a multiple of $4$. Letting
$\theta_i=2i\pi/N_\theta$, $i=0,\ldots, N_\theta-1$ and denoting
$\Delta\theta=2\pi/N_\theta$, we use the compound trapezoid formula
to get $\bar{B}(\cdot,\cdot)$, i.e., for $v,w\in V$,
\begin{eqnarray*}
\tilde{B}(v,w)&=&-\int_0^{2\pi}(D^\alpha_\theta v,D^\alpha_{\theta+\pi} w)\tilde{M}(\theta)d\theta+c(v,w)\\
&\approx&-\Delta\theta\sum\limits_{i=0}^{N_\theta-1}(D^\alpha_{\theta_i}
v,D^\alpha_{\theta_i+\pi} w)\tilde{M}(\theta_i)+c(v,w):= \bar{B}(v,w).
\end{eqnarray*}
The fact that  $\tilde{M}(\theta)=\tilde{M}(\theta+\pi)$ and $N_\theta$ is a multiple of 4
guarantees that $\bar{B}(v,w)$ is a symmetric bilinear form as
well, i.e., $\bar B(v,w)=\bar B(w,v)$. By Parseval's formula,
we have
\begin{eqnarray}
(D^\alpha_\theta v,D^\alpha_{\theta+\pi} w)_\Omega&=&((2\pi i\xi_1\cos \theta+2\pi i\xi_2\sin\theta)^{2 \alpha} \mathcal{F}{v},\overline{\mathcal{F}{w}})\nonumber\\
&\leq&C
||v||_{H^{\alpha}(\Omega)}||w||_{H^{\alpha}(\Omega)}
\end{eqnarray} and
\begin{eqnarray}\label{eqn tttt}
\frac{d}{d\theta}(D^\alpha_\theta v,D^\alpha_{\theta+\pi} w)_\Omega&=&2\alpha\left((-2\pi i\xi_1\sin \theta+2\pi i\xi_2\cos\theta)(2\pi i\xi_1\cos \theta+2\pi i\xi_2\sin\theta)^{2 \alpha-1} \mathcal{F}{v},\overline{\mathcal{F}{w}}\right)\nonumber\\
&\leq&C
||v||_{H^{\alpha}(\Omega)}||w||_{H^{\alpha}(\Omega)}.
\end{eqnarray}
By
the error formula for the compound trapezoid formula, it is easy
to verify that
\begin{equation}\label{eqn cond 1}
|\tilde{B}(v,w)-\bar B(v,w)|\leq C\Delta\theta
||v||_{H^{\alpha}(\Omega)}||w||_{H^{\alpha}(\Omega)},
\end{equation}
where  $C$ is a positive constant
independent of $\theta, v$ and $w$.  Combining (\ref{eqn cond 1}) with (\ref{r5}) and (\ref{r6}), we know for sufficiently small $\Delta\theta$,
$$\bar B(v,v)\gtrsim ||v||^2_{H^\alpha(\Omega)},\quad v\in H_0^\alpha(\Omega),$$
\begin{equation}\label{eqn aaaaaaa5}
\bar{B}(v,w)\lesssim
||v||_{H^\alpha(\Omega)}||w||_{H^\alpha(\Omega)}, \quad
v,w \in H_0^\alpha(\Omega).
\end{equation}
By Lax-Milgram theorem, (\ref{eqn finite element new r2}) has a unique solution. The first Strang
lemma (see \cite{Ciarlet1}) holds here, i.e.,
 \begin{eqnarray*}
 ||u-\bar u_h||_{H^\alpha(\Omega)}&\lesssim& C\inf\limits_{v\in V}\left\{||u-v||_{H^\alpha(\Omega)}+\sup\limits_{w\in V}\frac{|B(v,w)-\bar B(v,w)|}{||w||_{H^\alpha(\Omega)}}\right\}\\
 &\lesssim&C\inf\limits_{v\in V}\left\{||u-v||_{H^\alpha(\Omega)}+\Delta\theta ||v||_{H^\alpha(\Omega)}\right\}.
 \end{eqnarray*}

Finally, the finite element approximation of (\ref{eqn source new}) is unitedly presented as: find ${u}_h\in V$ such that
 \begin{equation}\label{eqn finite element new}
 B(u_h,v)=(f,v),\quad \forall v\in V,
 \end{equation}
where $B(v,w)=-\int_0^{2\pi}(D^\alpha_\theta v,D^\alpha_{\theta+\pi}
w)M(\theta)d\theta+c(v,w)$, $M(\theta)$ is equal to a discrete form $\sum_{k=1}^Lp_k\delta(\theta-\theta_k)$ such that $B(\cdot,\cdot)$ is a symmetric
bilinear form,
\begin{equation}\label{r5r6}
B(v,v)\gtrsim ||v||^2_{H^{\alpha}(\Omega)},B(v,w)\lesssim ||v||_{H^{\alpha}(\Omega)}||w||_{H^{\alpha}(\Omega)},\quad v,w\in H^{\alpha}_0(\Omega),
\end{equation}
and $\int_0^{2\pi}M(\theta)d\theta\lesssim1$.
Specially for the cases mentioned above, the finite element problem (\ref{eqn finite element new}) represents problem (\ref{eqn finite element new r1}) if ${M}(\theta)=\tilde{M}(\theta)=\sum_{k=1}^Lp_k\delta(\theta-\theta_k)$ and problem (\ref{eqn finite element new r2}) if $M(\theta)=\Delta\theta\sum_{i=0}^{N_\theta-1}\delta(\theta-\theta_i)\tilde{M}(\theta_i)$.

\section{Multigrid algorithm}\label{sec3}
In this section, for (\ref{eqn finite element new}), we shall present our V-cycle multigrid
algorithm and a general framework for our convergence analysis.

Take $f_h\in V$ such that $(f_h,v)=(f,v)$, $\forall v\in V$ and
define a linear operator $A:V\rightarrow V$ as follows:
\begin{equation}\label{eqn variational form 2}
(Av,w)=B(v,w), \quad \forall v,w\in V.
\end{equation}
The finite element approximation of system (\ref{eqn finite element new}) can be
restated as to find $u_h\in V$ such that
\begin{equation}\label{M eqn problem111}
Au_h=f_h.
\end{equation}
In the following, we shall use the operator equation (\ref{M eqn problem111}) to construct our multigrid
algorithm. Since $B(v,w)$ is a symmetric
bilinear form, we know, by (\ref{r5r6}), that
$A:V\rightarrow V$ is symmetric positive definite with respect to
$(\cdot,\cdot)$, i.e.,
$$(Av,w)=(v,Aw),\quad v,w\in V; \quad (Av,v)>0,\quad 0\neq v\in V.$$
Then bilinear form
$$(v,w)_A:=(Av,w),\quad  v,w\in V,$$
also induces an inner product on $V$. Set norm
$$||v||_A=(Av,v)^{1/2},\quad v\in V.$$
By (\ref{r5r6}), we have
\begin{equation}\label{eqn aaaa5}
||v||_A\approx ||v||_{H^\alpha(\Omega)}, \quad \forall v\in V.
\end{equation}
\subsection{Algorithm}\label{sec31}
Assume that the triangulation ${\cal T}_h$ of $\Omega$ is constructed by
a successive refinement process. To be precise,  let ${\cal T}_J ={\cal T}_h$ for
some $J>1$, and ${\cal T}_k$ for $ k\geq 0$ be a nested sequence of
quasi-uniform triangulations,  i.e., ${\cal T}_k=\{\tau_k^i\}$ consists of
simplexes $\tau_k^i$ of size $h_k$ such that $\Omega=\cup_i
\tau_k^i$; $\tau_{k-1}^l$ is a union of simplexes of $\tau_k^i$. We
further assume that there is a positive constant $\gamma<1$, independent of
$k$, such that $h_k$ is proportional to $\gamma^{k}$ and the simplexes in ${\cal T}_1$ are of diameter $\approx 1$.

For each partition ${\cal T}_k$, we may define finite element spaces
$V_k$ by
\begin{equation}
\label{multi space} V_k=\{v\in
C^0{(\bar\Omega)}:v|_{\partial\Omega}=0, v|_\tau\in
P_l(\tau),\forall \tau\in {\cal T}_k\}.
\end{equation}
Obviously, the following inclusion relation holds: $V_1\subset
V_2\subset\cdots\subset V_J=V.$
 Our V-cycle multigrid methods  are based on the subspace decomposition
 $V=V_1+V_2+\cdots+V_J.$

For each $k\in \{1,2,\ldots,J\}$, define projectors
$Q_k,P_k:V\rightarrow V_k$  by
$$(Q_kv,w)=(v,w),\quad (P_kv,w)_A=(v,w)_A,\quad v\in V, w\in V_k,$$
specially, set $Q_0:V\rightarrow V$ as $Q_0v=0$, and define the
linear operator $A_k:V_k\rightarrow V_k$
$$(A_kv,w)=(Av,w),\quad v,w\in V_k.$$
It is easy to verify that
\begin{equation}\label{M eqn important}
A_kP_k=Q_kA,\quad k=1,2,\ldots,J.
\end{equation}
It is obvious that $A_k$ is symmetric and positive definite with
respect to $(\cdot,\cdot)$. Denote $\lambda_k\in
\mathbb{R}$, $k=1,2,\ldots,J$, the maximal eigenvalue of $A_k$.

Let $u_k=P_ku_h$ and $f_k=Q_kf_h$, we may get the operator equation in subspace
\begin{equation}\label{eqn tt}
A_ku_k=f_k.
\end{equation}
Our multigrid algorithm is essentially an iterative
procedure in which the subspace equation (\ref{eqn tt}) is
approximately solved successively to get new approximations to
(\ref{M eqn problem111}) from old approximations. More precisely,
denote $R_k:V_k\rightarrow V_k$ the approximate inverse of $A_k$, and
$u^{old}$ the old approximation to $u$. Correcting the residual of
$u^{old}$ in $V_k$ gives
$$u^{new}=u^{old}+R_kQ_k(f_h-Au^{old}).$$
We take $R_k$ to be symmetric with respect to
$(\cdot,\cdot)$ such that
\begin{equation}\label{tem6}(R_kv,v)\approx \frac{1}{\lambda_k}(v,v),\quad \forall v\in V_k, k=1,2,\ldots,J.
\end{equation}
\begin{remark}\label{remark7}
In this paper, we have $h_1=O(1)$ and take $R_1=A^{-1}_1$.
By Lemma \ref{lemma V in H32}, (\ref{eqn aaaa5}) and the definition of norm $||\cdot||_{H^\mu(\Omega)}$, we know that $(v,v)\lesssim (A_1v,v)\lesssim h_1^{-2\alpha}(v,v)$, and  $\lambda_1=O(1)$. Then we have $(R^{-1}_1v,v)\approx\frac1\lambda_1(v,v).$
\end{remark}

Next we give our V-cycle multigrid algorithm. \\

{\bf V-cycle Multigrid Algorithm}.
 Let $u^0=0\in V$, assume that $u^k\in V$ has been obtained. Then $u^{k+1}$ is generated by
 \begin{equation}\label{iteration}
 u^{k+1}=u^{k}+B_J(f_h- Au^k),
 \end{equation} where $B_J$ is defined inductively: Let $B_1=A_1^{-1}$, and assume
 that $B_{k-1}:V_{k-1}\rightarrow V_{k-1}$ has been defined; then for $g\in V_k$, $B_k: V_{k}\rightarrow V_{k}$ is defined as follows:\\
\indent Step 1. $v^1=R_kg;$\\
\indent Step 2. $v^2=v^1+ B_{k-1}Q_{k-1}(g-A_kv^1);$ \\
\indent Step 3. $B_kg=v^2+R_k(g-A_kv^2).$

 \subsection{A general framework}
  For the V-cycle multigrid method, we have
 $$u_h-u^{k+1}=(I-B_JA)(u_h-u^{k}).$$
Denote
  \begin{equation}\label{eqn EJ}
 E_J=(I-T_J)(I-T_{J-1})\cdots (I-T_1),\quad E^*_J=(I-T_1)\cdots(I-T_{J-1})(I-T_J)
 \end{equation}
  with $T_1=P_1, T_k=R_kA_kP_k, k=2,3,\ldots,J$.
Then we have $(I-B_JA)=E_J E^*_J$. Define the operator norm as
$$||E_J||_A=\sup\limits_{v\in V}\frac{||E_Jv||_A}{||v||_A}.$$
 It is easy to see that  $E_J^*$
is the $(\cdot,\cdot)_A$-adjoint of $E_J$, i.e.,
$$(E_Jv,w)_A=(v,E_J^*w)_A, \quad v,w\in V$$
and that
$$||E_J||_A=||E^*_J||_A,\quad ||E_JE_J^*||_A\leq ||E_J||^2_A.$$

The main work in this paper is to establish the contraction
property: there is a constant $0<\delta<1$ independent of the mesh size and mesh level
such that
\begin{equation}\label{eqn rate}
||E_J||_A\leq \sqrt{\delta}.
\end{equation}
By (\ref{eqn rate}),  we may obtain $||u_h-u^k||_A\leq
\delta^{k}||u_h-u^0||.$

\begin{remark}\label{remark 6}
For the V-cycle multigrid method, the spectral radius of the iterative matrix
$\rho=\rho(I-B_JA)\leq \delta$. It is known that the condition
number $\kappa(B_JA)\leq \frac{1+\rho}{1-\rho}\leq
\frac{1+\delta}{1-\delta}$ and $B_JA$
is self-adjoint and positive with respect to inner product
$(\cdot,\cdot)_A$. The $\delta$'s
independence of the mesh size implies that $B_J$ is a good
preconditioner for $A$ which can be used to design efficient
preconditioned conjugate gradient methods.
\end{remark}

Define $K_0$ and $K_1$ as two smallest positive constants satisfying
the following conditions:

1. For any $v\in V$, there exists a decomposition
$v=\sum_{i=1}^Jv_i$ for $v_i\in V_i$ such that
\begin{equation}\label{eqn condition for K0}
\sum\limits_{i=1}^J(R_i^{-1}v_i,v_i)\leq K_0(Av,v).
\end{equation}

2. For any $S\subset \{1,2,\ldots,J\}\times \{1,2,\ldots,J\}$ and
$v_i,w_i\in V$ for $i=1,2,\ldots,J$,
\begin{equation}\label{eqn condition for K1}
\sum\limits_{(i,j)\in S}(T_iv_i,T_jw_j)_A\leq
K_1\left(\sum\limits_{i=1}^J(T_iv_i,v_i)_A\right)^{\frac12}\left(\sum\limits_{j=1}^J(T_jw_j,w_j)_A\right)^{\frac12}.
\end{equation}
The estimate of the upper bound of $||E_J||_A$ relies on the
following lemma:
\begin{lemma}\cite{Bramble1, Xu1}\label{lemma  for convergence} Let $E_J$ be defined by (\ref{eqn EJ}). We have
$$||E_J||_A\leq 1-\frac{2-\omega_1}{K_0(1+K_1)^2},$$
where $\omega_1=\max\limits_k\rho(R_kA_k)$, $\rho(R_kA_k)$ denotes
the spectral radius of $R_kA_k$.
\end{lemma}

The estimate of the parameter $\omega_1$ is straightforward. Since $R_1=A^{-1}_1$, $\rho(R_1A_1)=1$. From
(\ref{tem6}), for $v\in V_k$ ($k=2,\ldots,J$)
$$\frac{C_1}{\lambda_k}(v,v)\leq (R_kv,v)\leq\frac{C_2}{\lambda_k}(v,v),$$
and furthermore
 \begin{equation}\label{eqn ttt}
 (R_kA_kv,v)_A=(R_kA_kv,A_kv)\leq\frac{C_2}{\lambda_k}(A_kv,A_kv)\leq{C_2}(v,A_kv)=(v,v)_{A},
 \end{equation}
 where the last inequality is obtained from that $A_k$ is symmetric positive matrix and $\lambda_k$ is the maximal eigenvalue of $A_k$.
Combining (\ref{eqn ttt}) with the fact that $R_kA_k$ is symmetric
with respect to inner product $(\cdot,\cdot)_A$, we have
$\rho(R_kA_k)\leq C_2$. Taking $R_k$ such that $C_2$ is suitably small can guarantee
the $\omega_1< 2.$

Next, we shall estimate the parameters $K_1,K_2$. The
following Lemma is helpful for the analysis.
\begin{lemma}\cite{Bramble1, Xu1}\label{lemma Cauchy-Schwarz}
Let $\epsilon=(\epsilon_{ij})\in R^{J\times J}$ be a nonnegative
symmetric matrix, with components $\epsilon_{ij}$ being the smallest
constant satisfying
\begin{equation}\label{eqn Cauchy Schwarz}
(T_iv,T_jw)_A\leq \epsilon_{ij}
(T_iv,v)_A^{\frac12}(T_jw,w)_A^{\frac12},\quad \forall v,w\in V.
\end{equation}
Then we have
$$
K_1\leq \rho(\epsilon),
$$
where $\rho(\epsilon)$ denotes the spectral radius of matrix
$\epsilon$. Furthermore, if $\epsilon_{ij}\lesssim  \gamma^{|i-j|}$ for
some $\gamma\in (0,1)$, then $\rho(\epsilon)\lesssim
(1-\gamma)^{-1}$.
\end{lemma}

\section{Convergence Analysis}\label{sec4}
We here first introduce two interpolation norms and relevant Sobolev spaces (see e.g., \cite{Tartar1}).
Let $\Lambda$ be a
domain in $\mathbb{R}^2$. For integer $m$, denote by $||\cdot||_{\tilde{H}^{m}(\Lambda)}$ the Sobolev norm of integer order $m$, i.e.,
$$||v||_{\tilde{H}^m(\Lambda)}:=\left(\sum\limits_{|l|\leq m}||D^lv||_{L^2(\Lambda)}^2\right)^{1/2},$$
with $l=(l_1,l_2)$, $|l|=l_1+l_2$ and
$D^l=(\frac{\partial}{\partial x})^{l_1}(\frac{\partial}{\partial y})^{l_2}$. Let $\mu>0$ be a non-integer and $0<s<1$,  $n$ is a non-negative integer such that $n<\mu<n+1$. We introduce the interpolation norms
\begin{equation}\label{norm 1}||v||_{\tilde{H}^{\mu}(\Lambda)}:=\left(\int_0^{\infty}\tilde{K}(v,t) t^{-2\mu-1}dt\right)^{1/2},||v||_{\hat{H}^{s}(\Lambda)}:=\left(\int_0^{\infty}\hat{K}(v,t) t^{-2s-1}dt\right)^{1/2}
\end{equation}
where
$$\tilde{K}(v,t):=\inf_{w\in \tilde{H}^{n+1}(\Lambda)}\left(||v-w||^2_{\tilde{H}^n(\Lambda)}+t^2||w||^2_{\tilde{H}^{n+1}(\Lambda)}\right),$$
$$\hat{K}(v,t):=\inf_{w\in \tilde{H}_0^{1}(\Lambda)}\left(||v-w||^2_{L^2(\Lambda)}+t^2||w||^2_{\tilde{H}^{1}(\Lambda)}\right).$$
 Relevant Sobolev spaces are
\begin{equation}\label{eqn aaaaaa3}
{\tilde{H}}^{\mu}(\Lambda):=\{v\in
L^2(\Lambda);||v||_{{\tilde{H}}^{\mu}(\Lambda)} <\infty\},\quad{\hat H}^{s}(\Lambda):=\{v\in
L^2(\Lambda);||v||_{{\hat{H}}^{s}(\Lambda)} <\infty\}.
\end{equation}
Let $\Lambda_1$, $\Lambda_2$ be two domains in $\mathbb{R}^2$ with $\Lambda_1\subset \Lambda_2$, and then
\begin{eqnarray}\label{rr2}&&\left(\int_0^{\infty}\inf_{w\in \tilde{H}^{n+1}(\Lambda_1)}\left(||v-w||^2_{\tilde{H}^{n}(\Lambda_1)}+t^2||w||^2_{\tilde{H}^{n+1}(\Lambda_1)}\right) t^{-2\mu-1}dt\right)^{1/2}\nonumber\\
&\leq& \left(\int_0^{\infty}\inf_{w\in \tilde{H}^{n+1}(\Lambda_2)}\left(||(v-w)|_{\Lambda_1}||^2_{\tilde{H}^{n}(\Lambda_1)}+t^2||\,w|_{\Lambda_1}||^2_{\tilde{H}^{n+1}(\Lambda_1)}\right) t^{-2\mu-1}dt\right)^{1/2}\nonumber\\
&\leq& \left(\int_0^{\infty}\inf_{w\in \tilde{H}^{n+1}(\Lambda_2)}\left(||v-w||^2_{\tilde{H}^{n}(\Lambda_2)}+t^2||w||^2_{\tilde{H}^{n+1}(\Lambda_2)}\right) t^{-2\mu-1}dt\right)^{1/2}.
\end{eqnarray}
So we have, for $v\in \tilde{H}^{\mu}(\Lambda_2)$,
\begin{equation}\label{eqn r3}
||v||_{\tilde{H}^{\mu}(\Lambda_1)}\leq ||v||_{\tilde{H}^{\mu}(\Lambda_2)}.
\end{equation}
\begin{remark}\label{remark r1}
The following space relations can be found in literature: (1) $\mu>0$, ${\tilde{H}}^{\mu}(\mathbb{R}^2)$ and ${\tilde{H}}^{1}_0(\Omega)$ coincide with  ${H}^{\mu}(\mathbb{R}^2)$  and ${H}^{1}_0(\Omega)$ respectively;
(2) for $1/2<\mu<1$, ${\tilde{H}_0^\mu(\Omega)}$ coincides with $\hat{H}^{\mu}(\Omega)$ (see \cite{Lions1,Tartar1}); for $1/2<\mu<1$, ${\tilde{H}_0^\mu(\Omega)}$ coincides with  ${H}_0^{\mu}(\Omega)$ (this can been shown by (1), (2) and the definitions of the interpolation spaces).
\end{remark}
Combining with remark \ref{remark r1} and the well known interpolation property (see e.g., Lemma 22.3 in \cite{Tartar1}), we know, for $1/2<\mu\leq1$,
\begin{equation}\label{eqn aaaaaa2}
||(I-Q_k)v||\lesssim  h_k^\mu
||v||_{H^\mu(\Omega)}, \quad v\in H^\mu_0(\Omega).
\end{equation}

Now, we develop some results for the finite element
spaces $V_k,k\geq 1$.
 Let $\Omega'\subset \mathbb{R}^2$ be a
suitable polygonal domain such that $\Omega\subset\Omega'$ and
$\hbox{dist}(\partial \Omega',\Omega)>C$ for a positive $C$. ${\cal T}'_k$,
$k\geq 1$, are the quasi-uniform triangulations obtained by
extending ${\cal T}_k$ from $\Omega$ to $\Omega'$, that is, ${\cal T}'_k$ in
$\Omega$ coincides with ${\cal T}_k$. Furthermore we still make sure that
${\cal T}'_k=\{\tau_k^i\}$ consists of simplexes $\tau_k^i$ of size $h_k$.
Let $V'_k=\{v\in C^0{(\bar\Omega')}:v|_{\partial\Omega'}=0,
v|_\tau\in P_l(\tau),\forall \tau\in {\cal T}'_k\}.$ In the following,
for  $v\in V_k$, $v$ always denotes its extensions (on
$\Omega'$ and on $\mathbb{R}^2$), which is extended by zero outside
$\Omega$, and so we also have $v\in V'_k$.

\begin{lemma}\label{lemma r2}
Let $\mu>0$, $v\in {\tilde{H}}^{\mu}(\Omega')$ with $\hbox{supp}(v)\subset \Omega$ ($v$ also denotes its extension on $\mathbb{R}^2$ which is extended by zero outside
$\Omega'$). Then we have $||v||_{\tilde{H}^{\mu}(\Omega')}\approx ||v||_{H^{\mu}(\mathbb{R}^2)}.$
\end{lemma}
\proof For $\mu$ being a integer, the conclusion is direct. For the case that $\mu$ is not a integer, denote $n$ as a non-negative integer such that $n<\mu<n+1$.
 From (\ref{rr2}), $||v||_{\tilde{H}^{\mu}(\Omega')}\leq ||v||_{\tilde{H}^{\mu}(\mathbb{R}^2)}\approx||v||_{H^{\mu}(\mathbb{R}^2)}.$ Now we prove the converse relation. Let $\Lambda$ be a domain in $\mathbb{R}^2$ with $C^{n+1}-$smooth boundary such that $\Omega\subset\subset\Lambda\subset \Omega'$. Then by (\ref{eqn r3}), $v\in {\tilde{H}}^{\mu}(\Lambda)$.
Following the proof for the strong extension of Sobolev space (see e.g., Theorem 4.26 in \cite{Adams1}), we can show that there is a linear operator $E$ continuous from $\tilde{H}^{j}(\Lambda)$ into $\tilde{H}^{j}(\mathbb{R}^2)$ for integers $0\leq j\leq n+1$, such that $E(v|_{\Lambda})=v$. Then we have
\begin{eqnarray}\label{rr3}||v||_{\tilde{H}^{\mu}(\mathbb{R}^2)}&=&\left(\int_0^{\infty}\inf_{w\in \tilde{H}^{n+1}(\mathbb{R}^2)}\left(||v-w||^2_{\tilde{H}^{n}(\mathbb{R}^2)}+t^2||w||^2_{\tilde{H}^{n+1}(\mathbb{R}^2)}\right) t^{-2\mu-1}dt\right)^{1/2}\nonumber\\
&\leq& \left(\int_0^{\infty}\inf_{w\in \tilde{H}^{n+1}(\Lambda)}\left(||E(v|_{\Lambda}-w)||^2_{\tilde{H}^{n}(\mathbb{R}^2)}+t^2||Ew||^2_{\tilde{H}^{n+1}(\mathbb{R}^2)}\right) t^{-2\mu-1}dt\right)^{1/2}\nonumber\\
&\lesssim& \left(\int_0^{\infty}\inf_{w\in \tilde{H}^{n+1}(\Lambda)}\left(||v-w||^2_{\tilde{H}^{n}(\Lambda)}+t^2||w||^2_{\tilde{H}^{n+1}(\Lambda)}\right) t^{-2\mu-1}dt\right)^{1/2}\nonumber\\
&=&||v||_{\tilde{H}^{\mu}(\Lambda)},
\end{eqnarray}
where the last inequality is by the continuity of $E$. Combining with (\ref{eqn r3}), we obtain $||v||_{H^{\mu}(\mathbb{R}^2)}\approx||v||_{\tilde{H}^{\mu}(\mathbb{R}^2)}\lesssim ||v||_{\tilde{H}^{\mu}(\Omega')}.$ \eproof

\begin{lemma}\label{lemma V in H32}
For $0<\mu<3/2$, $v\in V_k$, we have
\begin{equation}\label{eqn a4}
||v||_{{{H}}^{\mu}(\mathbb{R}^2)}\lesssim  h_k^{-\mu}||v||,
\end{equation}
and then $V_k\subset {{H}}^{\mu}(\mathbb{R}^2)$.
\end{lemma}
\proof By $v\in V_k$, we known $v\in V'_k$,  $||v||_{\tilde{H}^{\mu}(\Omega')}\lesssim  h_k^{-\mu}||v||$ from \cite{Brammble,Xu1,Xu2} and further  $||v||_{H^{\mu}(\mathbb{R}^2)}\lesssim  h_k^{-\mu}||v||$ by Lemma \ref{lemma r2}. \eproof

Let $\beta$ be a positive with $\alpha+\beta<3/2$ and
$\alpha-\beta\geq 0$ in the rest of this paper. We have the
following results:
\begin{lemma}\label{lemma for condition 2}
It holds that
$$\quad (v,w)_A\lesssim  ||v||_{H^{\alpha+\beta}(\mathbb{R}^2)}||w||_{H^{\alpha-\beta}(\mathbb{R}^2)},\quad v,w\in V.$$
\end{lemma}
\proof Since $v,w\in V$, by Lemma \ref{lemma V in H32}, we know that $v,w\in H^{\alpha+\beta}({\mathbb{R}^2})$. Then
\begin{eqnarray*}
(v,w)_A&=&(Av,w)=B(v,w)\\
&=&-\int_0^{2\pi}\left(D^\alpha_\theta v,D^\alpha_{\theta+\pi} w\right)M(\theta)d\theta+c(v,w)\\
&=&-\int_0^{2\pi}\left(D^{\alpha+\beta}_\theta v,D^{\alpha-\beta}_{\theta+\pi} w\right)M(\theta)d\theta+c(v,w)\\
&\leq&\int_0^{2\pi}||D^{\alpha+\beta}_\theta v||_{L^2(\Omega)} ||D^{\alpha-\beta}_{\theta+\pi} w||_{L^2(\Omega)} M(\theta)d\theta+c||v||_{L^2(\Omega)}\,||w||_{L^2(\Omega)}\\
&\lesssim&||v||_{H^{\alpha+\beta}(\mathbb{R}^2)} || w||_{H^{\alpha-\beta}(\mathbb{R}^2)}+c||v||_{L^2(\Omega)}\,||w||_{L^2(\Omega)}\\
&\lesssim &||v||_{H^{\alpha+\beta}(\mathbb{R}^2)} ||
w||_{H^{\alpha-\beta}(\mathbb{R}^2)},\end{eqnarray*} where the third
equality is by Lemma \ref{lemma a7}, and the second inequality is by
Lemma \ref{lemma a9} and $\int_0^{2\pi}M(\theta)d\theta\lesssim1$. \eproof
\begin{lemma}\label{lemma Cauchy-Schwarz pp}
Let $i\leq j$, then
\begin{equation}\label{eqn Cauchy Schwarz pp}
(v,w)_A\lesssim   \gamma^{(j-i)\beta}
h_i^{-\alpha}h_j^{-\alpha}||v||\, ||w||,\quad v\in V_i,w\in V_j.
\end{equation}
Here we recall that $\gamma\in (0,1)$ is a constant such that
$h_k=O(\gamma^k)$.
\end{lemma}
\proof  For $v\in V_i,w\in V_j$, we know that
 \begin{eqnarray*}
 (v,w)_A&\lesssim & ||v||_{H^{\alpha+\beta}(\mathbb{R}^2)}||w||_{H^{\alpha-\beta}(\mathbb{R}^2)}
 \lesssim h_i^{-(\alpha+\beta)} ||v||\,h_j^{-(\alpha-\beta)}||w||\\
 &=&(h_j/h_i)^{\beta}h_i^{-\alpha}h_j^{-\alpha} ||v||\,||w||
 \lesssim  \gamma^{(j-i)\beta} h_i^{-\alpha}h_j^{-\alpha} ||v||\,||w||,
\end{eqnarray*}
where the first inequality is by Lemma \ref{lemma for condition 2},
the second inequality is by Lemma \ref{lemma V in H32}, and the last
inequality is by the relation $h_k\approx O(\gamma^k)$. \eproof

\begin{lemma}\label{lemma pp}
Let $W_i=(Q_i-Q_{i-1})V$, then
\begin{equation}\label{eqn pp}
(v,w)_A\lesssim \gamma^{|j-i|\beta}||v||_A||w||_A,  \quad \forall
u\in W_i,v\in W_j.
\end{equation}
\end{lemma}
\proof By (\ref{eqn aaaaaa2}) and (\ref{eqn aaaa5}), we have
$$||v||\lesssim h^\alpha_k||v||_A, \quad \forall v\in W_k.$$
Combining the above inequality with Lemma \ref{lemma Cauchy-Schwarz pp} gives the
lemma.

\begin{lemma}\label{lemma for Cauchy-Schwarz}
It holds that
\begin{equation}\label{eqn for Cauchy Schwarz 2}
(T_iv,T_jw)_A \lesssim \gamma^{|i-j|\beta}
(T_iv,v)_A^{\frac12}(T_jw,w)_A^{\frac12},\quad \forall v,w\in V.
\end{equation}
\end{lemma}
\proof It suffices to prove (\ref{eqn for Cauchy Schwarz 2}) holds
for $i\leq j$. Assume that $i\leq j$, and then for $v, w\in V$,
\begin{eqnarray}\label{eqn lq}
(T_iv,T_jw)_A&=&(R_iA_iP_iv,R_jA_jP_jw)_A\nonumber\\
&\lesssim &\gamma^{(j-i)\beta}
h_i^{-\alpha}h_j^{-\alpha}||R_iA_iP_iv||\, ||R_jA_jP_jw||,
\end{eqnarray}
where the inequality is by Lemma \ref{lemma Cauchy-Schwarz pp}.
$$||R_iA_iP_iv||^2=(R_iA_iP_iv,R_iA_iP_iv)\approx \frac{1}{\lambda_i}(R_iA_iP_iv,A_iP_iv)= \frac{1}{\lambda_i}(T_iv,v)_A,$$
where the second equality is by (\ref{tem6}) and the symmetry of
$R_k$. Then we obtain
\begin{equation}\label{eqn qq2}
||R_iA_iP_iv||\lesssim \lambda^{-1/2}_i(T_iv,v)^{1/2}_A,
\end{equation}
and similarly
\begin{equation}\label{eqn qq1}
||R_jA_jP_jw||\lesssim \lambda^{-1/2}_j(T_jw,w)^{1/2}_A.
\end{equation}
 For $v\in V_k$, we have
\begin{equation}\label{eqn aaaaaa5}
(Av,v)=||v||^2_A\approx ||v||^2_{H^\alpha(\Omega)}\lesssim
{h_k^{-2\alpha}}||v||^2,
\end{equation}
where the second equality is by (\ref{eqn aaaa5}) and the last inequality
is by Lemma \ref{lemma V in H32}. For $w\in V$, $v=(Q_k-Q_{k-1}) w\in
V_k$, by  (\ref{eqn aaaaaa2}), we have
\begin{equation}\label{eqn aaaaaa6}{h_k^{-2\alpha}}||v||^2\lesssim ||v||^2_{H^\alpha(\Omega)}\approx (Av,v).
\end{equation}
 By (\ref{eqn aaaaaa5}) and (\ref{eqn aaaaaa6}), it is not
hard to see that
\begin{equation}\label{range of eigenvalue 1}
\lambda_k\approx  {h_k^{-2\alpha}}, \quad  k=1,2,\ldots,J.
\end{equation}
Combining  (\ref{eqn lq}) with  (\ref{eqn qq2}), (\ref{eqn qq1}) and
(\ref{range of eigenvalue 1})  gives
\begin{equation*}
(T_iv,T_jw)_A \lesssim \gamma^{(j-i)\beta}
(T_iv,v)_A^{\frac12}(T_jw,w)_A^{\frac12},\quad \forall v,w\in V.
\end{equation*}
The Lemma is proved. \eproof
\begin{lemma}\label{lemma norm equivalence}
Let
\begin{equation}\label{eqn for norm }
||v||_M^2:=\sum\limits_{k=1}^J||(Q_k-Q_{k-1})v||^2_{A},
\end{equation}
and then for $v\in V$, we have
$$||v||_M\approx ||v||_{A}.$$
\end{lemma}
\proof It is not hard to see that the space $H^\alpha_0(\Omega)$ coincides with $\tilde{H}^\alpha(\Omega)$ in \cite{Oswald1}. Combining with Theorem 1 of \cite{Oswald1},  we know that $||w||^2_{H^{\alpha}(\Omega)}\approx\sum_{k=1}^\infty{h_k^{-2\alpha}}||(Q_k-Q_{k-1})w||^2$ holds for $w\in H^\alpha_0(\Omega)$.  For $v\in V$, $||(Q_k-Q_{k-1})v||^2\approx {h_k^{2\alpha}}||(Q_k-Q_{k-1})v||^2_{H^\alpha(\Omega)}$ by (\ref{eqn aaaaaa5}) and (\ref{eqn aaaaaa6}).  Combining with (\ref{eqn
aaaa5}) gives the lemma. \eproof
 \begin{theorem}\label{lemma estimate for K0 K1} We have
$$K_0\lesssim 1,\quad K_1\lesssim  1.$$
That is to say, our V-cycle multigrid method  is optimal, which means that the
convergence rate is independent of the mesh size and mesh level.
\end{theorem}
\proof For $v\in V$, decompose $v$ as $v=\sum_{i=1}^Jv_i$
with $v_i=(Q_i-Q_{i-1})v$. By (\ref{eqn aaaaaa2}) and (\ref{eqn
aaaa5}) we have $||v_i||\lesssim  h_i^{\alpha}||v_i||_A.$ Furthermore
combining (\ref{tem6}) with (\ref{range of eigenvalue 1}), we
have $(R_i^{-1}v_i,v_i)\lesssim  ||v_i||^2_A.$ Using Lemma
\ref{lemma norm equivalence} gives $K_0\lesssim 1$.  Finally Combining Lemma
\ref{lemma for Cauchy-Schwarz} with Lemma \ref{lemma Cauchy-Schwarz}
gives that $K_1\lesssim  1$. \eproof

\section{Implementation }
 For simplicity, in this section, we only consider $l=1$ in (\ref{multi space}), i.e., $V_k$,
$k=1,2,\ldots,J$, are the spaces consisting of the piecewise linear
polynomials. Let $\phi^i_k$, $i=1,\ldots,N_k$, be the nodal basis of the finite element space $V_k$.  The implementation  are a classical procedure in literature (see e.g., \cite{Bramble1}), and  we here only illustrate how to generate the stiff matrices of the finite element systems and how to choose  $R_k:V_k\rightarrow V_k$, $k=2,\ldots,J$, the
approximations of $A_k$.

\subsection{The stiffness matrices and $R_k$}
For $A_k$, denote its corresponding
stiffness matrix by $\tilde A_k\in \mathbb{R}^{N_k\times N_k}$ with
entries
\begin{equation}\label{tem1}
(\tilde A_k)_{ij}=B(\phi^i_k,\phi^j_k).
\end{equation}
Since $M$ has the discrete form $M(\theta)=\sum_{l=1}^L
p_l\delta(\theta-\theta_l)$,
$$(\tilde{A}_k)_{ij}=-\sum\limits_{l=1}^{L}p_l(D^{2\alpha-1}_{\theta_l} \phi^i_k,D_{\theta_l+\pi} \phi^j_k)+c(\phi^i_k,\phi^j_k).$$
We need only  discuss how to numerically compute
\begin{eqnarray}\label{tem2}I_\theta&=&(D^{2\alpha-1}_{\theta} \phi^i_k,D_{\theta+\pi} \phi^j_k)=(D^{-\nu}_{\theta} D_{\theta}\phi^i_k,D_{\theta+\pi} \phi^j_k)\nonumber \\
&=&\int_{\hbox{ssup}(\phi_k^j)}D^{-\nu}_{\theta}
D_{\theta}\phi^i_k\times D_{\theta+\pi} \phi^j_k dxdy
\end{eqnarray}
for a fixed $\theta$, where $\nu=(2-2\alpha)$, and then the
entries of the stiff matrices can be numerically computed. If
$\alpha=1$, the computation of the stiffness matrices  is easy, since
the original problem is an integer order one. Now we focus on the case
of $1/2<\alpha<1$. Define the index set $K_i$ as
$$K_i=\{l;\tau_k^l\in {\cal T}_k, \tau_k^l\subset \hbox{supp}(\phi^i_k)\}.$$
Then
\begin{eqnarray*}\label{tem3}
I_\theta&=&\sum\limits_{l\in K_j} \int_{\tau_k^l} D^{-\nu}_{\theta} D_{\theta}\phi^i_k\times D_{\theta+\pi} \phi^j_k dxdy\\
&=&\sum\limits_{l\in K_j}\sum\limits_{l'\in K_i}  \int_{\tau_k^l}
D^{-\nu}_{\theta} (\chi_{\tau_k^{l'}}D_{\theta}\phi^i_k)\times
D_{\theta+\pi} \phi^j_k dxdy,
\end{eqnarray*}
where for a set $S$
in $\mathbb{R}^2$,
$$\chi_S(x,y)=\left\{
\begin{array}{ll} 1, & if \ (x,y)\in S; \\
                  0, & otherwise. \end{array}\right. $$
Noting that $D_{\theta}(\phi^i_k)|_{\tau_k^{l'}}$,
$D_{\theta+\pi}(\phi^j_k)|_{\tau_k^{l}}$ are both constants, we numerically compute
\begin{equation}\label{tem4}
\int_{\tau_k^l}
D^{-\nu}_{\theta}\chi_{\tau_k^{l'}}(x,y)\times\chi_{\tau_k^{l}}(x,y)dxdy,
\end{equation}
and then $I_\theta$ can be computed.
\begin{figure}
\begin{center}
  \includegraphics[width=4cm]{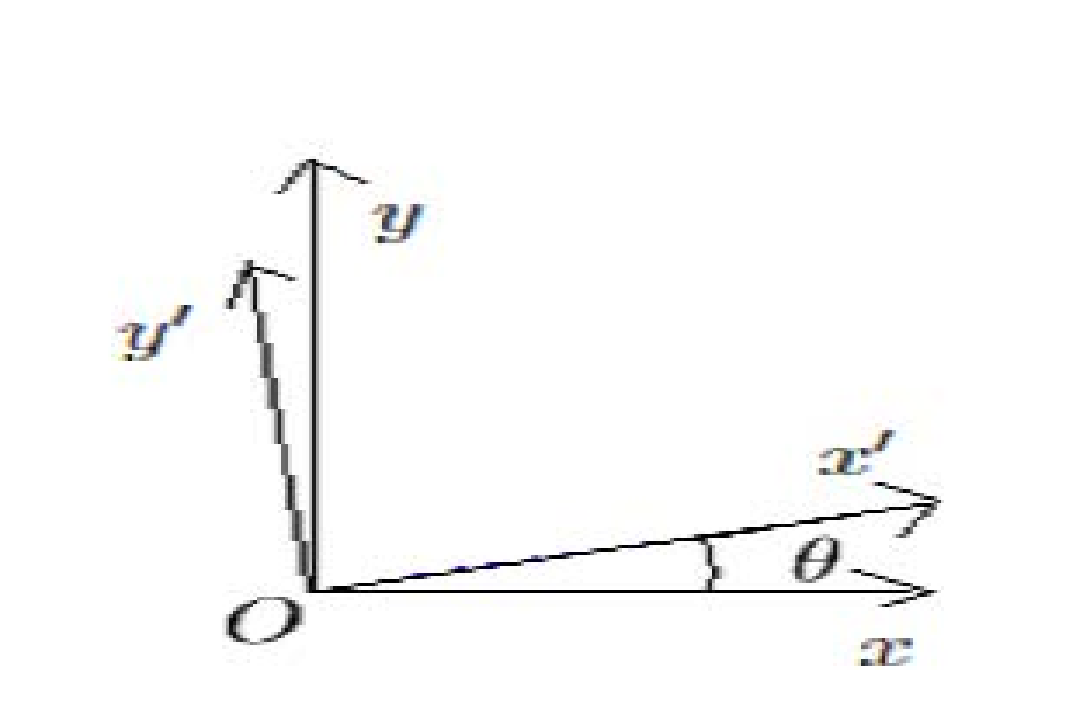}\hspace{1cm}\includegraphics[width=6cm]{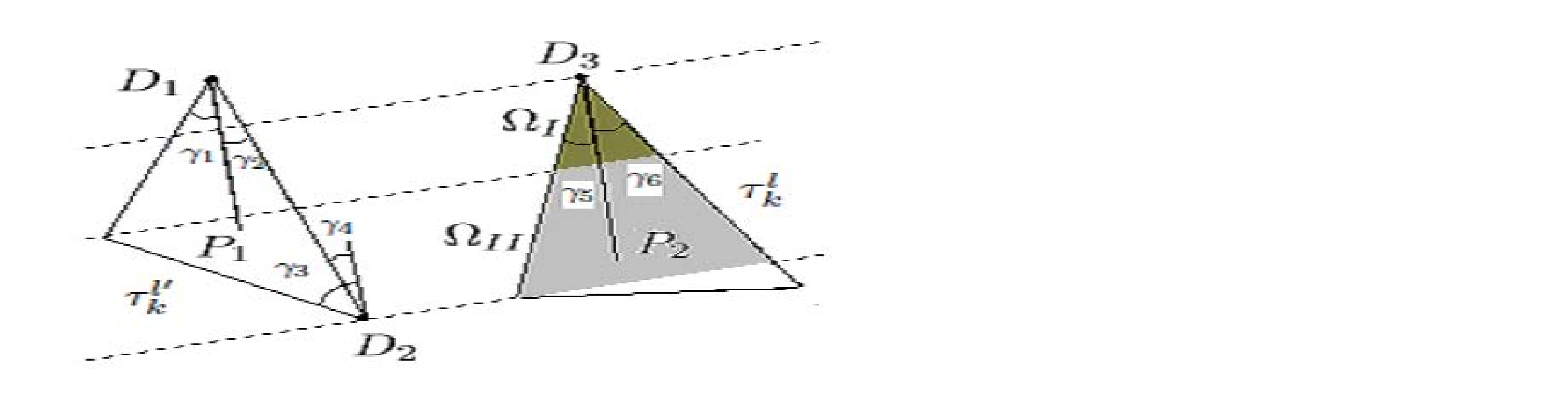}\\
 \caption{Illustration for computing $I_\theta$}\label{graph3}
 \end{center}
\end{figure}
Next we illustrate how to compute the integral in (\ref{tem4}) by an example.
On  the left of Figure \ref{graph3} is Cartesian coordinate systems
$xOy$ and $x'Oy'$, and the angle between axes $Ox$ and $Ox'$ is
$\theta$. On  the right of Figure \ref{graph3}, the two triangles
are $\tau_k^{l'}$ and $\tau_k^{l}$; $D_1,D_2,D_3$ denote the
corresponding vertices of the triangles; $\Omega_I,\Omega_{II}$
denote the corresponding shadow areas respectively; lines $D_1P_1$
and $D_3P_2$ are both Parallel to axis $Oy'$; $\gamma_{1},
\gamma_{2}, \gamma_{3}, \gamma_{4},  \gamma_{5}, \gamma_{6}$ are
correspondent angles. Denote the coordinates of $D_1,D_2$ and $D_3$
under coordinate system $x'Oy'$ by $(x'_{11},y'_{12})$,
 $(x'_{21},y'_{22})$ and $(x'_{31},y'_{32})$ respectively.
 Then we have
 \begin{eqnarray*}\label{tem5}
 &&\int_{\tau_k^l} D^{-\nu}_{\theta}\chi_{\tau_k^{l'}}(x,y)\times\chi_{\tau_k^{l}}(x,y)dxdy\\
 &=&\int_{\Omega_I} D^{-\nu}_{\theta}\chi_{\tau_k^{l'}}(x,y)dxdy+\int_{\Omega_{II}} D^{-\nu}_{\theta}\chi_{\tau_k^{l'}}(x,y)dxdy\\
 &=&\int_{\Omega_I} {_{-\infty}} D^{-\nu}_{x'}\chi_{\tau_k^{l'}}(x',y')dx'dy'+\int_{\Omega_{II}} {_{-\infty}}D^{-\nu}_{x'}\chi_{\tau_k^{l'}}(x',y')dx'dy'\\
 &=&\int_{\Omega_I} \frac{1}{\Gamma(\nu+1)}\left(x'-x'_1+(y'_1-y')\tan\gamma_1\right)^\nu dx'dy'\\
& &\quad -\int_{\Omega_I} \frac{1}{\Gamma(\nu+1)}\left(x'-x'_1-(y'_1-y')\tan\gamma_2\right)^\nu dx'dy'\\
& &\int_{\Omega_{II}} \frac{1}{\Gamma(\nu+1)}\left(x'-x'_2+(y'-y'_2)\tan\gamma_3\right)^\nu dx'dy'\\
& &\quad -\int_{\Omega_{II}}
\frac{1}{\Gamma(\nu+1)}\left(x'-x'_2+(y'-y'_2)\tan\gamma_4\right)^\nu
dx'dy'.
 \end{eqnarray*}
The last four integrals above can be computed directly. Finally
we know that the entries of the
stiffness matrices can be numerically computed.

We choose $R_k$ as
\begin{equation}\label{tem7}
R_kv=\frac{1}{\tilde{\lambda}_k}\sum\limits_{i=1}^{N_k}(v,\phi_k^i)\phi_k^i,\quad
v\in V_k,
\end{equation}
with $\tilde{\lambda}_k\approx \lambda_kh_k^2$. Define mass matrix
$M_k\in \mathbb{R}^{N_k\times N_k}$ with entries
$$(M_k)_{ij}=(\phi_k^i,\phi_k^j).$$
For $v\in V_k$, denote $\tilde{v}\in \mathbb{R}^{N_k}$ the vector of
coefficients of $v$ in the basis $\{\phi^i_k\}_{i=1}^{N_k}$.  It is
 known that $\tilde{v}^T M_k\tilde{v}\approx h_k^2\tilde{v}^T\tilde{v}$ and $\tilde{v}^T
M_k^2\tilde{v}\approx h_k^2\tilde{v}^TM_k\tilde{v}$. Hence we have
\begin{equation}\label{eqn new}(R_kv,v)\approx \frac{1}{\tilde{\lambda}_k}\tilde{v}^TM_k^2\tilde{v}\approx \frac{1}{{\lambda}_k}\tilde{v}^TM_k\tilde{v}\approx \frac{1}{\lambda_k}(v,v),
\end{equation}
which  means (\ref{tem6}) holds.
In the numerical tests, we take $\widetilde{\lambda}_k=\frac{3}2(\tilde A_{k})_{ii}$, $k=2,\ldots,J.$
It is not hard to verify that $
(\tilde A_k)_{ii}\approx h_k^{2-2\alpha}\approx h_k^2\lambda_k.
$

\subsection{Computation complexity}
For the numerical approximation of SFPDEs, one of the key issues
is how to reduce the computation complexity. We confine ourself
to the case that $\Omega$ is a square domain, and of course the technique
here is also helpful for effectively designing schemes for the case
that $\Omega$ is a general domain.

\begin{figure}
\begin{center}
 \includegraphics[width=8cm]{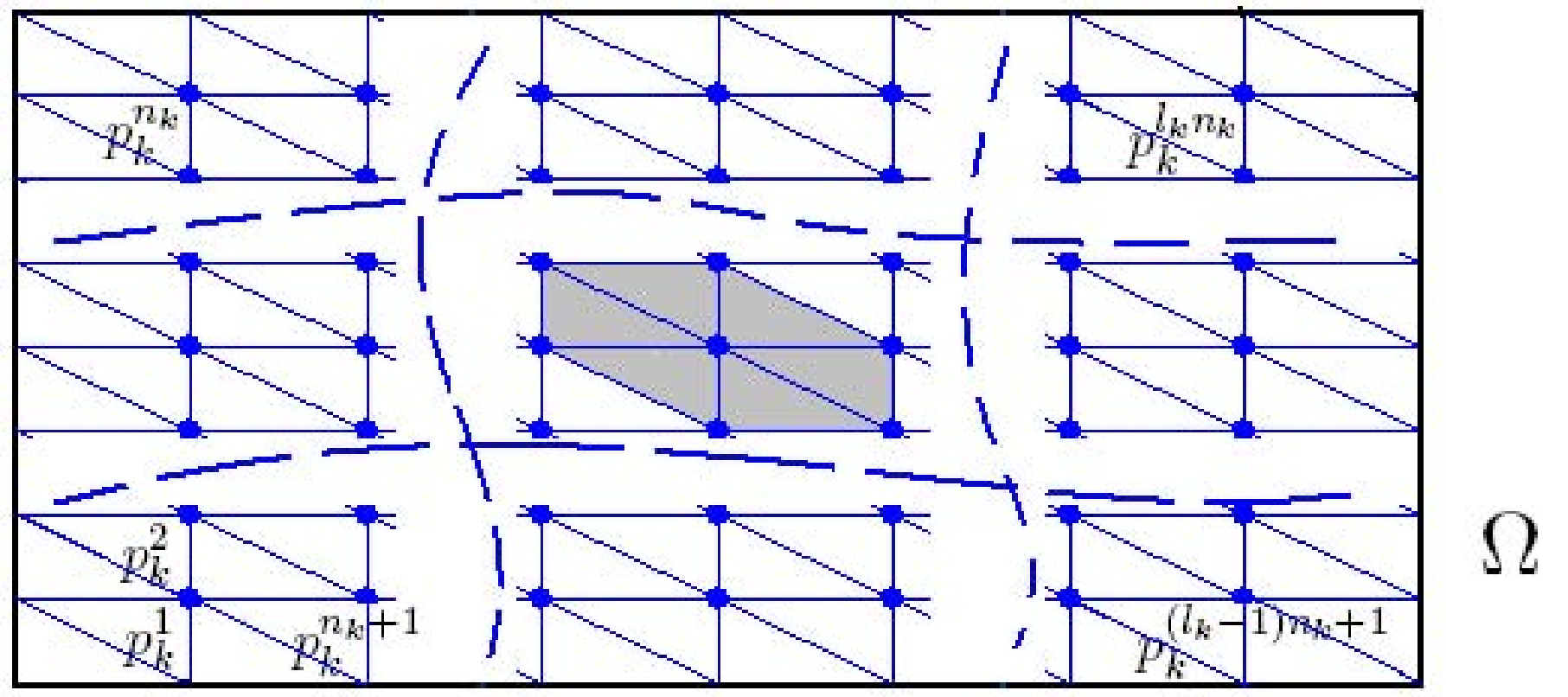}\\
 \caption{Uniform triangulation}\label{graph1}
 \end{center}
\end{figure}

The triangulations ${\cal T}_k$, $k=1,2,\ldots,J$ are those in Figure
\ref{graph1}, where dashed curve denote the ellipsis, $n_k=n_02^{k}-1,l_k=l_02^k-1$ with positive
integers $n_0,l_0$, and $p_k^{m}$, $m=1,\ldots,n_kl_k$
are the interior points. The finite element space $V_k=\{v\in
H^1_0(\Omega): v|_\tau\in P_1(\tau),\forall \tau\in {\cal T}_k\}$. Let
$\phi_k^{m}=\phi_k^{m}(x,y)$, $m=1,\ldots,n_kl_k$, be the nodal
basis functions, i.e., $\phi_k^{m}$ is a piecewise linear polynomial
whose values are 1 at $p_k^{m}$ and zeros at other nodes (including
interior and exterior nodes).

Denote
$U=(U_{1},U_{2}\ldots,U_{n_k},\ldots,U_{2n_k,},\ldots,U_{l_kn_k})^T$.
Next we discuss how to effectively conduct the
multiplication of matrix $\tilde{A}_k$ and  vector $U$.
Let $\nu=(\nu_1,\nu_2,\ldots,\nu_{(2n_k-1)l_k-n_k+1})^T\in
\mathbb{R}^{(2n_k-1)l_k-n_k+1}$ with
$$\nu_{(2n_k-1)j+i}=B(\phi_k^{1},\phi_k^{jn_k+i}),\ i=1,\ldots,n_k, j=0,\ldots,l_k-1,$$
$$\nu_{(2n_k-1)j-i+2}=B(\phi_k^{i},\phi_k^{jn_k+1}),\ i=2,\ldots,n_k, j=1,\ldots,l_k-1.$$
Define a symmetric Toeplitz matrix
$$\tilde{\tilde{A}}=\left(
              \begin{array}{ccccc}
                \nu_1 & \nu_2 & \cdots & \nu_{(2n_k-1)l_k-n_k} & \nu_{(2n_k-1)l_k-n_k+1} \\
                \nu_2 & \nu_1 & \cdots & \nu_{(2n_k-1)l_k-n_k-1} & \nu_{(2n_k-1)l_k-n_k} \\
                \vdots  & \vdots  & \ddots & \vdots              & \vdots \\
                \nu_{(2n_k-1)l_k-n_k}  & \nu_{(2n_k-1)l_k-n_k-1} & \ldots & \nu_1 & \nu_2 \\
                \nu_{(2n_k-1)l_k-n_k+1} & \nu_{(2n_k-1)l_k-n_k} & \ldots & \nu_2 & \nu_1 \\
              \end{array}
            \right).
$$
 Toeplitz matrix, also called diagonal-constant matrix, is a matrix in which each descending diagonal from left to right is a constant.

For any  $i,j$ with $1\leq i\leq j\leq n_kl_k$, let
$d_i,r_i,d_j,r_j$ be nonnegative integers satisfying $i=n_kd_i+r_i$,
$j=n_kd_j+r_j$, $1\leq r_i,r_j\leq n_k$. Let
$j'=(d_j-d_i),i'=|r_j-r_i|$, and then by the property of the
operator $B(\cdot,\cdot)$, it is easy to see that
$$
B(\phi_k^{i},\phi_k^{j})=\left\{
\begin{array}{ll}
  B(\phi_k^1,\phi_k^{j'n_k+i'+1})=\nu_{j'(2n_k-1)+i'+1}, & if\ r_j\geq r_i; \\
  B(\phi_k^{i'+1},\phi_k^{j'n_k+1})=\nu_{j'(2n_k-1)-i'+1}, & if\ r_j< r_i.
\end{array}\right.
$$
And thereby any component of matrix $\tilde{\tilde{A}}$ is also one of
vector $\nu$. Define sets
$$\mathcal{I}_m=\{m(2n_k-1)+1,m(2n_k-1)+2,\ldots,m(2n_k-1)+n_k\},\quad m=0,1,\ldots,l_k-1$$
and $\mathcal{I}=\bigcup\limits_{0\leq m\leq l_k-1} \mathcal{I}_m$.
We have the relation
\begin{equation}\label{Nu eqn 4}
\tilde{A}_k=\tilde{\tilde{A}}_{\mathcal{I},\mathcal{I}},
\end{equation}
where $\tilde{\tilde{A}}_{\mathcal{I},\mathcal{I}}$ denotes the
sub-matrix of  $\tilde{\tilde{A}}$ which consists of entries
$\tilde{\tilde{A}}_{ij}$ of $\tilde{\tilde{A}}$ indexed by $i,j\in
\mathcal{I}$. Denote $U'\in \mathbb{R}^{(2n_k-1)l_k-n_k+1}$ as
$$U'=(U_{1},\ldots,U_{n_k},\overbrace{0,\ldots,0}\limits^{n_k-1},U_{n_k+1},\ldots,U_{2n_k},\overbrace{0,\ldots,0}\limits^{n_k-1},U_{2n_k+1},\ldots,U_{l_kn_k}).$$
It is not hard to see that
$$\tilde{A}_kU=(\tilde{\tilde{A}}U')_{\mathcal{I}},$$
where for a given vector $v$, $v_{\mathcal{I}}$ denotes the vector which
consists of entries $v_i$ indexed by $i\in \mathcal{I}$.  So the
multiplication of the matrix $\tilde{A}_k$ and any vector $U\in
R^{n_kl_k}$ can be obtained by conducting the multiplication of the
Toeplitz matrix $\tilde{\tilde A}$ and $U'\in R^{(2n_k-1)l_k-n_k+1}.$
The multiplication of a Toeplitz matrix in $\mathbb{R}^{n\times n}$
and a vector in $\mathbb{R}^n$ can be done with computation complexity
$O(n\log n)$. Recall that $N_J=n_Jl_J$ denotes the number of the
unknowns in the finite element problem (\ref{M eqn problem111}), and
then by the above analysis, we conclude that for the V-cycle
multigrid methods developed in Section \ref{sec4}, each iteration
needs computation complexity $O(N_J\log N_J)$.

\subsection{Numerical results}
In this section, we shall present some numerical results to confirm our theoretical findings.  In our
numerical test, we take $n_0=l_0=4$, and take $N_\theta=4(n_J+1)$ if
$M$ is a continuous function.

We shall check our V-cycle multigrid method and the
preconditioned conjugate gradient algorithm (PCG) with $B_J$ as the
preconditioner. Meanwhile, the numerical result for the conjugate gradient algorithm (CG) is also presented
for comparison. Our tests are carried out using Matlab software.
The stopping criterion of the algorithm is
 $$||u^k-u^{k-1}||_\infty\leq 10^{-6}.$$

We  present two examples: one is with the probability measure $\tilde{M}$
having a discrete form and the other with $\tilde{M}$ being a continuous
function. Table \ref{table 1} and Table
\ref{table 2} list
 the numerical results for Example \ref{exam 1} and Example \ref{exam 2} respectively, where "DOFs" denotes the degree of freedoms
and "Iter" denotes the iterative steps on each level.
 It is seen that the numbers of iterations of our V-cycle multigrid and PCG  per level
are  bounded independent of the mesh size and mesh level, which confirms our theoretical results.

\begin{example}\label{exam 1} Let $\Omega=[0,2]\times[0,2]$, the equation to be solved is
\begin{equation}\label{problem exam 1}
-\frac14( _{-\infty}D^{1.5}_{x}+  {_{x}}D^{1.5}_{\infty}+
{_{-\infty}D^{1.5}_{y}}+{_{y}D^{1.5}_{\infty}})u=1.
\end{equation}

\end{example}

\begin{table}[!hbp]
\begin{center}
\begin{tabular}{|c|c|c|c|c|c|}
\hline
{Level} & \multirow{2}{*}{DOFs} & {V-cycle} & {PCG} & {CG} \\
\cline{3-5}
            J           &            & Iter  & Iter      & Iter  \\
\hline
               4        &  4096      & 13     & 7          & 58    \\
               5        &  16384     & 13     & 6         & 72    \\
               6        &  65536     & 13     & 7         & 118   \\
               7        &  262144    & 13     & 7        & 197   \\
               8        &  1048576   & 12     & 7         & 313   \\
\hline
\end{tabular}
\caption{Numerical results for Example \ref{exam 1}.}
\label{table 1}
\end{center}
\end{table}

\begin{figure}
\begin{center}
 \includegraphics[width=6cm]{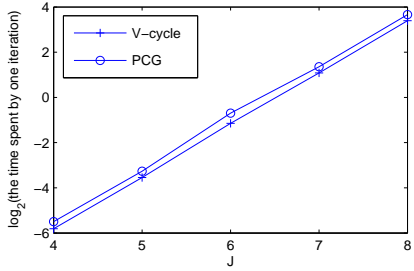}\\
 \caption{The CPU time per iteration}\label{figure times}
 \end{center}
\end{figure}

\begin{example}\label{exam 2} Let $\Omega=[0,2]\times[0,2]$ and $\tilde{M}(\theta)=1$. The equation to be solved is
\begin{equation}\label{problem exam 2}
-D^{1.5}_{\tilde{M}} u=1.
\end{equation}
\end{example}

 We choose smooth $f(x,y)$ in the examples such that the solutions have singularity near the boundaries. The computation complexity of our multigrid methods are shown in figure \ref{figure times}, where "Time" denotes
the CPU time (in seconds)  spent by one iteration.
As can be seen from the  figure \ref{figure times}, the CPU time of each iteration
ia almost linear with respect to the degree of freedoms.
So the computation complexity of our multigrid method is also optimal.

\begin{table}[!hbp]
\begin{center}
\begin{tabular}{|c|c|c|c|c|c|}
\hline
\multirow{2}{*}{Level} & \multirow{2}{*}{DOFs} & V-cycle & PCG  & CG \\
\cline{3-5}
                       &    & Iter & Iter  & Iter  \\
\hline
               4        &  4096     & 11  & 6  & 50 \\
               5        &  16384    & 11  & 6   & 69 \\
               6        &  65536    & 11  & 6   & 115 \\
               7        &  262144   & 11  & 6  & 184 \\
               8        &  1048576  & 11  & 6 & 307 \\
\hline
\end{tabular}
\caption{Numerical results for Example \ref{exam 2}.}
\label{table 2}
\end{center}
\end{table}


\renewcommand{\thesection}{Appendix }
\section{}
\renewcommand{\thesection}{A}
The Fourier analysis plays critical roles in this paper: for $g\in
L^1(\mathbb{R}^2)$, the Fourier transform of $g$ is the function
$\mathcal{F}g$  defined on (the dual of) $\mathbb{R}^2$
by
$$\mathcal{F}g(\xi_1,\xi_2)=\int_{\mathbb{R}^2}e^{-2i\pi(x\xi_1+y\xi_2)}g(x,y)dxdy,$$
where $i$ denotes the imaginary unit; for $g\in L^2(\mathbb{R}^2)$,
the Fourier transform $\mathcal{F}g$ of $g$ is defined in
the following distribution sense (see \cite{Tartar1}):
$$(\mathcal{F}g,v)=(g,\mathcal{F}v), \quad\forall v\in C^\infty_0(\mathbb{R}^2),$$
and more precisely, $\mathcal{F}$ is an isometry from
$L^2(\mathbb{R}^2)$ into itself, which satisfies Parseval's formula (see
\cite{Rudin1})
$$||\mathcal{F}g||=||g||$$
and
$$(v,\overline{w})=(\mathcal{F}v,\overline{\mathcal{F}w}),$$
where $\overline{z}$ denotes the complex conjugate of the complex
number $z$. The Fourier transform of the $\mu$th order fractional
derivative consists of the complex in the form $(i\kappa)^\mu$ with
$\mu>0, \kappa\in \mathbb{R}$ (see \cite{Podlubny1}). So it may be a
multi-valued function. To guarantee the Fourier transform to be
univalent, we express complex variable  $z=|z|\exp({i\theta})$,
$-\pi\leq\theta<\pi$, where
$\exp({i\theta})=\cos\theta+i\sin\theta$, $|z|$ and $\theta$
respectively denote the modulus and the argument of $z$. Then
$$(i\kappa)^\mu=(\hbox{sign}(\kappa)i|\kappa|)^\mu=(|\kappa|\exp({i\hbox{sign}(\kappa)\pi/2}))^\mu=|\kappa|^\mu \exp({i\mu\hbox{sign}(\kappa)\pi/2}),$$
$$(-i\kappa)^\mu=|\kappa|^\mu \exp({-i\mu\hbox{sign}(\kappa)\pi/2}).$$
It is easy to see that, for $\mu>0$,
\begin{equation}\label{eqn aaaa6}
\overline{(-i\kappa)^\mu}=(i\kappa)^\mu,\quad \forall \kappa\in
\mathbb{R}.
\end{equation}


\begin{thebibliography}{bib}
\small

\bibitem{Adams1} R. A. Adams, Sovolev Spaces, Academic Press, New York, 1975.

\bibitem{Bramble1}
J. Bramble,  Multigrid Methods, Pitman, Boston, 1993.

\bibitem{Brammble} J. Bramble, J. Pasciak and J. Xu, The analysis of multigrid algorithm with nonested spaces or noninherited quadratic forms, Math. Comp., 56 (1991), pp. 1-34.

\bibitem{Brammble1} J. Bramble, J. Pasciak and P. Vassilevski, Computational scales of Sovolev norms with application to preconditioning, Math. Comp., (69) 1999, pp. 463-480.

\bibitem{Beumer1} B. Beumer, M. Kov¨¢cs and M. M. Meerschaert, Numerical solutions for fractional reaction diffusion equations, Comput. Math. Appl., 55 (2008), pp. 2212-2226.

\bibitem{Susanne1}S. Brenner and L. Scott, The
Mathematical Theory of Finite Element Methods, Springer-Verlag, New
York, 1994.

\bibitem{Butzer1} P.  Butzer and H. Berens, Semi-groups of Operators and Approximation, Springer-Verlag,
Berlin and New York, 1967.

\bibitem{Chan1} R. Chan and X. Jin, An Introduction to Iterative Toeplitz Solvers, SIAM, Philadelphia, 2007.
\bibitem{Chan2} R. Chan and M. Ng, Conjugate gradient methods for Toeplitz systems, SIAM Rev., 38 (1996), pp. 427-482.

\bibitem{Ciarlet1} P. Ciarlet,  The Finite Element Methods for Elliptic Problems,   North-Holland, New York,
1978.


\bibitem{Cui1} M. Cui, Compact finite difference method for the fractional diffusion equation, J. Comput. Phys., 228 (2009), pp. 7792-7804.

\bibitem{Deng1} W. Deng, Finite element method for the space and time fractional Fokker-Planck equation, SIAM J. Numer. Anal., 47 (2008), pp. 204-226.

\bibitem{Dryja1}  M. Dryja and O. Widlund,  Towards a unified theory of domain decomposition algorithms for elliptic problems, in Domain decomposition Method for Partial Differential Equations, Society for Industrial and Applied Mathematics, Philadelphia, PA, 1990. Written in English.

\bibitem{Ervin0} V. Ervin, N. Heuer and J. Roop, Numerical approximation of a time dependent, nonlinear, space-fractional diffusion equation, SIAM J. Numer. Anal., 45
(2007), pp. 572-591.

\bibitem{Ervin1} V.  Ervin and J. Roop, Variational formulation for the
stationary fractional advection dispersion equation, Numer. Meth.
P.D.E.,  22 (2006), pp. 558-576.



\bibitem{Ervin2} V. Ervin and J. Roop, Variational solution of the fractional advection dispersion equation on bounded domains in $R^d$ ,
Numer. Meth. P.D.E., 23 (2007), pp. 256-281.

\bibitem{Lei1} S.  Lei and H.  Sun, A circulant preconditioner for fractional diffusion equations, J. Comput. Phys., 242 (2013), pp. 715-725.

\bibitem{Langlands1} T.  Langlands and B. Henry, The accuracy and stability of an implicit solution method for the fractional diffusion equation, J. Comput. Phys., 205
(2005), pp. 719-736.

\bibitem{Li1} X. Li and C. Xu, The existence and uniqueness of the weak solution of the space-time fractional diffusion equation and a spectral method approximation,
Commun. Comput. Phys., 8 (2010), pp. 1016-1051.

\bibitem{Lin1} F.  Lin, S. Yang and X.  Jin,  Preconditioned iterative methods for fractional diffusion equation, J. Comput. Phys., 256 (2014), pp. 109-117.

\bibitem{Lin2} Y. Lin and C. Xu, Finite difference/spectral approximations for the time-fractional diffusion equation, J. Comput. Phys., 225 (2007), pp. 1533-1552.

\bibitem{Lions1} J. L. Lions, E. Magenes, Nonhomogeneous Boundary Value Problems and Applications I, Springer, Berlin (1972)



\bibitem{Liu1} F. Liu, V. Anh and I. Turner, Numerical solution of the space fractional Fokker-Planck equation, J. Comput. Appl. Math., 166 (2004), pp. 209-219.

\bibitem{Meerschaert0} M. Meerschaert, D.  Benson and B. Baumer, Multidimensional advection and fractional dispersion, Phys. Rev. E., 59 (1999), pp. 5026-5028.

\bibitem{Meerschaert00} M.  Meerschaert, J. Mortensen and H. Scheffler, Vector Grunwald formula for fractional derivatives, Fract.
Calc. Appl. Anal., 7 (2004), pp. 61-82.

\bibitem{Meerschaert1} M.  Meerschaert and C. Tadjeran, Finite difference approximations for two-sided space-fractional partial differential equations, Appl. Numer. Math., 56 (2006), pp. 80-90.


\bibitem{Meerschaert2} M.  Meerschaert, H. Scheffler and C. Tadjeran, Finite difference methods for two-dimensional fractional dispersion equation, J. Comput. Phys., 211 (2006), pp. 249-261.

\bibitem{Metzler1} R. Metzler and J. Klafter, The random walk's guide to anomalous diffusion: A fractional dynamics approach, Phys. Rep., 339 (2000), pp. 1-77.

\bibitem{Nezza1} E. Di. Nezza, G. Palatucci and E. Valdinoci, Hitchhiker's guide to the fractional Sobolev spaces, Bull. Sci. Math., 136 (2012), pp. 521-573.

\bibitem{Oswald1} P. Oswald, Multilevel norms for $H^{-1/2}$, Computing, 61 (1998), pp. 235-255.

\bibitem{Pang1} H. Pang and H. Sun, Multigrid method for fractional diffusion equations, J. Comput. Phys., 231 (2012), pp. 693-703.

\bibitem{Podlubny1} I. Podlubny, Fractional Differential Equations, Academic Press, New York, 1999.

\bibitem{Roop1} J. Roop,  Computational aspects of FEM approximation of fractional advection
dispersion equations on bounded domains in R2, Journal of
Computational and Applied Mathematics, 193 (2006), pp. 243-268.

\bibitem{Rudin1} W. Rudin, Real and Complex Analysis, McGraw-Hill, New York, 1987.

\bibitem{Samko1}  S.  Samko, A.  Kilbas and O. Marichev, Fractional Integrals and Derivatives: Theory and
Applications, Gordon and Breach, New York, 1993.

\bibitem{Sousa1} E. Sousa, Finite difference approximates for a fractional advection diffusion problem, J. Comput. Phys. 228 (2009), pp. 4038-4054.


\bibitem{Servadei1} R. Servadei and E. Valdinoci, Variational methods for non-local operators of elliptic type,
Discrete Contin. Dyn. Syst., 33 (2013), pp. 2105-2137.

\bibitem{Smith1} B. Smith, P. Bjorstad and W. Gropp, Domain Decomposition: Parallel Multilevel Methods for Elliptic Partial Differential Equations, Cambridge University Press, 1996.

\bibitem{Tadjeran1} C. Tadjeran, M.  Meerschaert and H.  Scheffler, A second-order accurate numerical approximation for the fractional diffusion equation, J. Comput. Phys., 213 (2006), pp. 205-213.

 \bibitem{Tartar1}   L. Tartar, An Introduction to Sobolev Spaces and Interpolation Spaces. Lecture Notes of the Unione Matematica Italiana 3, Springer-Verlag, Berlin Heidelberg, 2007.


 \bibitem{Wang1}   H. Wang, K. Wang and T. Sircar, Adirect $O(N\log^2 N)$ finite difference method for fractional diffusion equations, J. Comput. Phys., 229 (2010), pp. 8095-8104.

\bibitem{Wang2} H. Wang and K. Wang, An $O(N\log^2 N)$ alternating-direction finite difference method for two-dimensional fractional diffusion equations, J. Comput. Phys., 21 (2011), pp 7830-7839.

\bibitem{Wang3} H. Wang and N. Du, Fast alternating-direction finite difference methods for three-dimensional space-fractional diffusion equations, J. Comput. Phys., 258 (2014), pp 305-318.

\bibitem{Wang4}  K. Wang and H. Wang, A fast characteristic finite difference method for fractional advection-diffusion equations, Adv. Water Resour, 34 (2011), pp. 810-816.

\bibitem{Wang5}  H. Wang and N. Du, A superfast-preconditioned iterative method for steady-state space-fractional diffusion equations, J. Comput. Phys., 240 (2013), pp. 49-57.

\bibitem{Wang6} H. Wang and N. Du, A fast finite difference method for three-dimensional time-dependent space-fractional diffusion equations and its efficient implementation, J. Comput. Phys., 253 (2013), pp. 50-63.

\bibitem{Xu1} J.  Xu, Iterative methods by space decomposition and subspace correction, SIAM Rev., 34 (1992), pp. 581-613.

\bibitem{Xu2} J. Xu, Theory of multilevel methods, Ph.D. thesis, Cornell University, Ithaca, NY, Rep, AM-48, Pennsylvania State University,
University Park, PA, 1989.

\bibitem{Zhou1} Z.  Zhou and H. Wu, Finite element multigrid method for the boundary value problem of fractional advection dispersion equation, J. Appl. Math., Volume 2013, Article ID 385463, 8 pages.


\end{thebibliography}
\end{document}